# THE GROWTH OF ADDITIVE PROCESSES

By Ming Yang

*Columbia University*

Let $X_t$ be any additive process in $\mathbb{R}^d$. There are finite indices $\delta_i$, $\beta_i$, $i=1,2$ and a function $u$, all of which are defined in terms of the characteristics of $X_t$, such that

$$\liminf_{t \to 0} u(t)^{-1/\eta} X_t^* = \begin{cases} 0, & \text{if } \eta > \delta_1, \\ \infty, & \text{if } \eta < \delta_2, \end{cases}$$

$$\limsup_{t \to 0} u(t)^{-1/\eta} X_t^* = \begin{cases} 0, & \text{if } \eta > \beta_2, \\ \infty, & \text{if } \eta < \beta_1, \end{cases} \quad \text{a.s.,}$$

where $X_t^* = \sup_{0 \le s \le t} |X_s|$. When $X_t$ is a Lévy process with $X_0 = 0$, $\delta_1 = \delta_2$, $\beta_1 = \beta_2$ and $u(t) = t$. This is a special case obtained by Pruitt. When $X_t$ is not a Lévy process, its characteristics are complicated functions of $t$. However, there are interesting conditions under which $u$ becomes sharp to achieve $\delta_1 = \delta_2$, $\beta_1 = \beta_2$.

**1. Introduction.** A process $X_t$ with independent increments, rcll (right-continous with left limits) paths and values in $\mathbb{R}^d$ is called additive if $X_t$ is continuous in probability and $X_0 = 0$. Additive processes represent a large family of nonhomogeneous processes and intersect the entirety of Feller processes at the class of Lévy processes. Pruitt [6] defined an index $\delta$ for each Lévy process $X_t$ with $X_0 = 0$ and showed that $X_t$ satisfies the Hölder conditions: $\liminf_{t \to 0} t^{-1/\eta} X_t^* = 0$ or $\infty$ a.s. according as $\eta > \delta$ or $\eta < \delta$, where $X_t^* = \sup_{0 \le s \le t} |X_s|$. Its lim sup analogue was obtained by Blumenthal and Getoor [1] with an index $\beta$. Both results have their additive process counterparts. We define in terms of the characteristics of an additive process $X_t$ a nondecreasing continuous function $u$ with $u(0) = 0$ and four finite indices $\delta_i, \beta_i, i = 1, 2$ such that

(1.1) $$\liminf_{t \to 0} u(t)^{-1/\eta} X_t^* = \begin{cases} 0, & \text{if } \eta > \delta_1, \\ \infty, & \text{if } \eta < \delta_2, \end{cases}$$

$$\limsup_{t \to 0} u(t)^{-1/\eta} X_t^* = \begin{cases} 0, & \text{if } \eta > \beta_2, \\ \infty, & \text{if } \eta < \beta_1, \end{cases} \quad \text{a.s.}$$

In the case of Lévy processes, $u(t) = t$, $\delta_1 = \delta_2 = \delta$, $\beta_1 = \beta_2 = \beta$. Schilling [7] studied form (1.1) with $u(t) = t$ for a class of Feller processes. The issue of defining $u$ other than the indices arises when $X_t$ is nonhomogeneous. We cannot define $u$ to be "$t$" or any particular function holding for all additive processes. For example, continuous maps $B_t : \mathbb{R}_+ \to \mathbb{R}^d$ are additive processes (deterministic) but it is obvious that $u(t) = B_t^* = \max_{0 \le s \le t} |B_s|$. Thus, $u$ depends on $X_t$. We can also define two finite indices $\delta, \beta$ and four functions $\overline{v}, \underline{v}, \overline{u}, \underline{u}$ (not necessarily monotone) in terms of the characteristics of an additive process $X_t$ such that with probability 1

(1.2a)
$$\lim_{t \to 0} t^{-1/\eta} X_{\overline{v}(t)}^* = \infty, \quad \text{if } \eta < \delta,$$
$$\liminf_{t \to 0} t^{-1/\eta} X_{\underline{v}(t)}^* = 0, \quad \text{if } \eta > \delta,$$

(1.2b)
$$\limsup_{t \to 0} t^{-1/\eta} X_{\overline{u}(t)}^* = \infty, \quad \text{if } \eta < \beta,$$
$$\lim_{t \to 0} t^{-1/\eta} X_{\underline{u}(t)}^* = 0, \quad \text{if } \eta > \beta.$$

In many cases $\overline{v}/\underline{v} \le 1, \overline{u}/\underline{u} \le 1$ hold automatically. Otherwise we can always define two functions $v(\eta, t), u(\eta, t)$ in terms of the characteristics of $X_t$ such that with probability 1

(1.2c)
$$\liminf_{t \to 0} t^{-1/\eta} X_{v(\eta', t)}^* = \begin{cases} 0, & \text{if } \eta \wedge \eta' > \delta, \\ \infty, & \text{if } \eta \vee \eta' < \delta, \end{cases}$$
$$\limsup_{t \to 0} t^{-1/\eta} X_{u(\eta', t)}^* = \begin{cases} 0, & \text{if } \eta \wedge \eta' > \beta, \\ \infty, & \text{if } \eta \vee \eta' < \beta. \end{cases}$$

Equation (1.2c) is a simple consequence implied by (1.2a), (1.2b). Are there functions $v_i, v_s$ (not necessarily monotone) and indices $\delta, \beta \in (0, \infty)$ such that

(1.3)
$$\liminf_{t \to 0} t^{-1/\eta} X_{v_i(t)}^* = \begin{cases} 0, & \text{if } \eta > \delta, \\ \infty, & \text{if } \eta < \delta, \end{cases}$$
$$\limsup_{t \to 0} t^{-1/\eta} X_{v_s(t)}^* = \begin{cases} 0, & \text{if } \eta > \beta, \\ \infty, & \text{if } \eta < \beta, \end{cases} \quad \text{a.s.?}$$

That is the question we are trying to get into. If $\delta_1 = \delta_2, \beta_1 = \beta_2$ in (1.1), (1.3) follows with $v_i = v_s = u^{-1}$ the inverse of $u$. If $\overline{v}/\underline{v} \le 1, \overline{u}/\underline{u} \le 1$, (1.3) holds for any functions $v_i, v_s$ satisfying $\overline{v} \le v_i \le \underline{v}, \overline{u} \le v_s \le \underline{u}$. Equation (1.3) is an accurate statement that increases the degree of technicality in defining desired quantities. Refer to the information in Section 5 for Schilling's work on (1.3).

This paper is organized as follows. Section 2 contains the background on additive processes and some technical results needed later on. In Section 3 we begin with the proof of (1.1) and then turn to the issue that $\delta_1 = \delta_2, \beta_1 = \beta_2$. In Section 4 we establish (1.2a), (1.2b) and find the cases in which $\overline{v}/\underline{v} \leq 1, \overline{u}/\underline{u} \leq 1$ hold. In Section 5 we show that $u$ in (1.1) can be represented as $Ee(X_t^*)$ for some bounded function $e$. ($e$ can be characterized as the benchmark function up to a log log term for the law of the iterated logarithm.) Finally, Section 6 leaves some existence questions in check toward the settlement of (1.3).

*Some terminology.* Two positive functions $f_1$ and $f_2$ are said to be *comparable*, written as $f_1 \approx f_2$, if $f_1/f_2$ is trapped inside a finite positive interval. A nondecreasing right-continuous function $\phi$ with $\phi(t) > 0, t > 0, \phi(0) = 0$ is called quasiconvex (resp. moderate) if there are two constants $\rho, \sigma \in (0, \infty)$ such that $\phi(t_2)/\phi(t_1) \geq \rho(t_2/t_1)^\sigma$ [resp. $\phi(t_2)/\phi(t_1) \leq \rho(t_2/t_1)^\sigma$] whenever $0 < t_1 < t_2$. The exponent $\sigma$ is not unique. In this paper the term *inverse* refers to the right-continuous inverse. $\phi$ is quasiconvex (moderate) if and only if its inverse is moderate (quasiconvex). Typically, $t^p(\log(1/t))^\kappa, t^p(\log\log(1/t))^\kappa, t^p(\log\log\log(1/t))^\kappa, p > 0, \kappa \in \mathbb{R}$, and so on, along with their inverses are both quasiconvex and moderate. $(\log(1/t))^{-\kappa}$, $(\log\log(1/t))^{-\kappa}, \kappa > 0$, and so on (their inverses) are moderate (quasiconvex) but, however, not quasiconvex (moderate). A function $c:(0,1) \to (0,1)$ is called slow if $\liminf_{r \to 0} c(r)r^{-\eta} > 0$ for all $\eta > 0$, equivalently $\lim_{t \to 0} t^\eta/c(t) = 0$ for all $\eta > 0$. Moderate functions $(\log(1/t))^{-p}, (\log\log(1/t))^{-p}, p > 0$, and so on, as well as constant functions are slow.

**2. Characteristics of additive processes.** Let $X_t$ be an additive process in $\mathbb{R}^d$. There are two measures and two kernels: (the jump measure) $\mu = \sum_{t \geq 0} 1(\triangle X_t \neq 0)\delta_{(t, \triangle X_t)}$ on $\mathbb{R}_+ \times \mathbb{R}^d$, where $\delta_a$ is the Dirac point mass at $a \in \mathbb{R}_+ \times \mathbb{R}^d$; (the intensity measure) $\nu(B) = E\mu(B), B \in \mathcal{B}(\mathbb{R}_+ \times \mathbb{R}^d)$; $\mu_t(A) = \mu([0,t] \times A) = \sum_{s \leq t} 1(\triangle X_s \in A, \triangle X_s \neq 0)$, $A \in \mathcal{B}(\mathbb{R}^d)$; $\nu_t(A) = \nu([0,t] \times A) = E\mu_t(A)$. $\nu_t$ is a Lévy measure for fixed $t$. If $A^c$ contains an open ball with center at $0$, $\nu_t(A)$ is a nondecreasing continuous function in $t$. Thus, $\nu_t$ is a nondecreasing continuous Lévy kernel. Conversely, any nondecreasing continuous Lévy kernel $\nu_t$ gives rise to a unique additive process $X_t$ up to an independent continuous additive process. The characteristic function for $X_t$ takes the form $E\exp\{i\langle \lambda, X_t\rangle\} = e^{\Psi_t(\lambda)}, \lambda \in \mathbb{R}^d$, where $\Psi_t(\lambda) = i\langle B_t, \lambda\rangle - 2^{-1}\langle \lambda, Q_t\lambda\rangle + \int[e^{i\langle\lambda,x\rangle} - 1 - i\langle\lambda,x\rangle 1(|x| \leq 1)]\nu_t(dx)$. $B_t = (B_t^{(1)}, B_t^{(2)}, \ldots, B_t^{(d)}) \in \mathbb{R}^d$ is continuous with $B_0 = 0$. $Q_t = (q_{ij}(t))_{d \times d}$ is a nonnegative definite symmetric $d \times d$ matrix, which defines a centered Gaussian process. For fixed $\lambda$, $\langle \lambda, Q_t\lambda\rangle$ is a nondecreasing continuous function in $t$ with $\langle \lambda, Q_0\lambda\rangle = 0$. Thus, the $C_t^{(i)} = q_{ii}(t)$ are nondecreasing continuous functions with $q_{ii}(0) = 0$. $q_{ij}(t), i \neq j$, the elements off the diagonal are continuous functions of bounded variation null at 0 because they are the predictable quadratic covariation processes of a $d$-dimensional continuous Gaussian martingale. The characteristics of the $i$th component $X_t^{(i)}$ of $X_t$ are $B_t^{(i)}, C_t^{(i)} = q_{ii}(t), \nu_t^{(i)}(B) = \nu_t(\{x \in \mathbb{R}^d : x_i \in B\}), B \in \mathcal{B}(\mathbb{R})$, respectively. Let $X_t$ be a real-valued additive process with $E \exp\{i\lambda X_t\} = e^{\Psi_t(\lambda)}, \lambda \in \mathbb{R}$ where $\Psi_t(\lambda) = i\lambda B_t - 2^{-1}\lambda^2 C_t + \int(e^{i\lambda x} - 1 - i\lambda x 1(|x| \leq 1))\nu_t(dx)$. Define for $r > 0, t \geq 0$

$$G_t(r) = \int_{|x|>r} \nu_t(dx),$$

(2.1)

$$K_t(r) = r^{-2}\left[C_t + \int_{|x|\leq r} x^2 \nu_t(dx)\right],$$

$$M_t(r) = r^{-1}\left|B_t + \int_{1<|x|\leq r\vee 1} x\nu_t(dx) - \int_{r\wedge 1<|x|\leq 1} x\nu_t(dx)\right|,$$

(2.2)

$$M_t^*(r) = \max_{0\leq s\leq t} M_s(r),$$

(2.3) $\quad y_t(r) = G_t(r) + K_t(r) + M_t^*(r).$

For any process $X_t$ in $\mathbb{R}^d$ with *additive components*, define

(2.4) $$y_t(r) = \sum_{i=1}^d y_t^{(i)}(r),$$

where the $y_t^{(i)}(r)$ are given by (2.3) for their respective components $X_t^{(i)}$ of $X_t$. Since each $X_t^{(i)}$ is continuous in probability, $y_t(r)$ is nondecreasing continuous in $t$ for each fixed $r > 0$ with $y_0(r) = 0$. While every additive process in $\mathbb{R}^d$ must have additive components, a process with additive components does not necessarily have independent increments. There are an infinite number of processes with additive components having identical marginals $(B_t^{(i)}, C_t^{(i)}, \nu_t^{(i)}), 1 \leq i \leq d$, some of which are additive in $\mathbb{R}^d$ including the one whose components are independent of one another. If $X_t$ is a Lévy process in $\mathbb{R}^d$, $y_t(r) = th(r)$ where $h$ is the same function as defined in [6]. $y_t(r)$ has a doubling property. That is, for all $\theta > 1, r > 0, t \geq 0$,

(2.5) $\quad (3\theta^2)^{-1} y_t(r) \leq y_t(\theta r) \leq 2y_t(r).$

The proof goes as follows: If $M_t(r)$ in (2.2) is nondecreasing in $t$, $M_t(r) = M_t^*(r)$ in which case by (2.3) of [6] in continuous time, for all $\theta > 1, r > 0, t \geq 0$, $(2\theta^2)^{-1} y_t(r) \leq y_t(\theta r) \leq 2y_t(r)$. In the matter of a few lines one covers the general case for arbitrary $M_t(r)$ with a left-side constant to decrease by one-sixth.

LEMMA 2.1. *Let $X_t$ be a process in $\mathbb{R}^d$ with additive components and $y_t(r)$ the function in (2.4). Then for all $r > 0, t \geq 0$,*

(2.6)
$$P(X_t^* \geq r) \leq \pi_d y_t(r), \qquad P(X_t^* \leq r) \leq A_k(d) y_t(r)^{-k/2},$$
$$k = 1, 2, \ldots,$$

*where $\pi_d = aK(d), a = 2^{-1}(3+\sqrt{5}), K(d) = 3d^2, d > 1, K(1) = 1$ and $A_k(d) = (18\sqrt{2dk})^k$.*

PROOF. The proof is essentially one dimensional and similar to that of (3.2) of [6]. Let $X_t$ be a real additive process with the Lévy–Itô decomposition $X_t = X_t^r + Y_t^r$ at the level $r$ where $Y_t^r$ is the step process constituted by only those jumps of $X_t$ with size bigger than $r$. The number of such jumps up to time $t$ follows a Poisson distribution with mean $G_t(r)$. Decompose $X_t^r$ further into an independent sum of two martingales, one continuous, one purely discontinuous, as $X_t^r = EX_t^r + X_t^{r,c} + X_t^{r,d}$. The quadratic variation information shows $E(X_t^{r,d})^2 = \int_{|x| \leq r} x^2 \nu_t(dx)$ and $E(X_t^{r,c})^2 = C_t$. Thus, $\operatorname{Var} X_t^r = C_t + \int_{|x|\leq r} x^2 \nu_t(dx) = r^2 K_t(r)$. Subtracting the exponents for $Y_t^r, X_t^{r,c}, X_t^{r,d}$ collectively from $\Psi_t(r)$ gives $EX_t^r = B_t + \int_{1 < |x| \leq r \vee 1} x \nu_t(dx) - \int_{r \wedge 1 < |x| \leq 1} x \nu_t(dx)$, or $|EX_t^r| = rM_t(r)$.

*The first inequality in* (2.6). Define $A = (Y_s^r \neq 0 \text{ for some } s \in (0, t])$, the event that there is at least one jump with size greater than $r$ up to time $t$. Then $P(A) = 1 - e^{-G_t(r)}$. Obviously, $A^c \cap (X_t^* \geq r) \subset (X_t^{r*} \geq r)$. It follows that

$$\begin{aligned}P(X_t^* \geq r) &= P((X_t^* \geq r) \cap A) + P((X_t^* \geq r) \cap A^c) \\ &\leq P(A) + P(X_t^{r*} \geq r) = 1 - e^{-G_t(r)} + P(X_t^{r*} \geq r) \\ &\leq G_t(r) + P(X_t^{r*} \geq r).\end{aligned}$$

By the continuous version of Kolmogorov's inequality (a special case of Doob's maximal inequality),

$$P\left(\sup_{0 \leq s \leq t} |X_s^r - EX_s^r| \geq (1 - a^{-1})r\right) \leq ar^{-2} \operatorname{Var} X_t^r = aK_t(r),$$

where $(1 - a^{-1})^{-2} = a = 2^{-1}(3+\sqrt{5})$. If $M_t^*(r) \geq a^{-1}$, $P(X_t^* \geq r) \leq 1 \leq aM_t^*(r) \leq ay_t(r)$. If $M_t^*(r) < a^{-1}$, $|EX_s^r| = rM_s(r) \leq rM_t^*(r) < a^{-1}r$ for all $s \in [0, t]$, which implies that $P(X_t^{r*} \geq r) \leq P(\sup_{0 \leq s \leq t}|X_s^r - EX_s^r| \geq (1-a^{-1})r) \leq aK_t(r)$. Thus, $P(X_t^* \geq r) \leq G_t(r) + aK_t(r) \leq ay_t(r) = \pi_1 y_t(r)$.

*The second inequality in* (2.6). Let $D_r = P(X_t^* \leq r)$. We show that $D_{2r} \leq 18\sqrt{2}y_t(r)^{-1/2}$ first. The concentration function for a real-valued r.v. $X$ is defined as $Q(X;r) = \sup_{x \in \mathbb{R}} P(x \leq X \leq x+r), r > 0$. Let $X$ be an infinitely divisible random variable having characteristic function $E\exp\{i\lambda X\} = \exp\{i\lambda b - 2^{-1}\lambda^2\sigma^2 + \int(e^{i\lambda x} - 1 - i\lambda x 1(|x| \leq 1))\nu(dx)\}, \lambda \in \mathbb{R}$ and define $q(r) = r^{-2}\sigma^2 + \int (x/r)^2 \wedge 1 \nu(dx), r > 0$. Then $Q(X;r) \leq \sqrt{2\pi}q(r)^{-1/2}$. This inequality can be found in [4], Chapter 15, page 408. Suppose that $K_t(2r) \geq 2(G_t(2r) + M_t^*(2r))$. Then $K_t(2r) \geq 2^{-1}K_t(2r) + G_t(2r) + M_t^*(2r) \geq 2^{-1}y_t(2r)$. Thus, $D_r \leq P(|X_t| \leq r) \leq Q(X_t; 2r) \leq \sqrt{2\pi}q(2r)^{-1/2} = \sqrt{2\pi}(G_t(2r) + K_t(2r))^{-1/2} \leq 2\sqrt{\pi}y_t(2r)^{-1/2}$. Suppose that $K_t(2r) \leq 2(G_t(2r) + M_t^*(2r))$. Consider the Lévy–Itô decomposition $X_t = X_t^{2r} + Y_t^{2r}$ at the level $2r$ as well as the three events: $A_1 = (Y_s^{2r} = 0, s \in (0,t]), A_2 = (X_t^{2r*} \leq r), A_3 = (X_t^* \leq r)$. Suppose that $A_1^c$ occurs and let $\tau \in (0,t]$ be the first jump time of $Y_s^{2r}$. Then $X_s = X_s^{2r}, s \in [0, \tau)$. If $\sup_{0 \leq s < \tau} |X_s| > r$, $A_3^c$ occurs. If $\sup_{0 \leq s < \tau} |X_s| \leq r$, $|X_{\tau-}| \leq r$. Therefore, $|X_\tau| = |X_{\tau-} + \triangle X_\tau| \geq |\triangle X_\tau| - |X_{\tau-}| > 2r - r = r$ and hence $A_3^c$ occurs again. Suppose that $A_2^c$ occurs and $A_1^c$ does not. Then $X_s^{2r} = X_s, s \in [0,t]$ and hence $A_2^c \setminus A_1^c \subset A_3^c$. We have shown that $A_1^c \cup A_2^c \subset A_3^c$, that is, $A_3 \subset A_1 \cap A_2$. Therefore, $D_r = P(A_3) \leq P(A_1) \wedge P(X_t^{2r*} \leq r)$. [In fact $A_3 = A_1 \cap A_2$ and $P(A_3) = P(A_1)P(A_2)$ but neither of them is needed in the proof.] If $G_t(2r) \geq cM_t^*(2r)$ for some number $c > 0$, $y_t(2r) \leq 3(1+c^{-1})G_t(2r)$. It follows that $D_r \leq P(A_1) = e^{-G_t(2r)} \leq (1+G_t(2r))^{-1} < G_t(2r)^{-1} \leq 3(1+c^{-1})y_t(2r)^{-1}$. If $G_t(2r) \leq cM_t^*(2r), y_t(2r) \leq 3(1+c)M_t^*(2r)$ and $K_t(2r) \leq 2(1+c)M_t^*(2r)$. If $M_t^*(2r) \leq a$ for some number $a > 1$, $D_r \leq 1 \leq aM_t^*(2r)^{-1} \leq 3(1+c)ay_t(2r)^{-1}$. If $M_t^*(2r) \geq a$, $|EX_{t^*}^{2r}| = 2rM_{t^*}(2r) \geq 2ra$ for some $t^* \in [0,t]$ satisfying $M_{t^*}(2r) = M_t^*(2r)$. [Note that $M_t(2r)$ is continuous in $t$.] Thus,

$$D_r \leq P(X_{t^*}^* \leq r) \leq P(X_{t^*}^{2r*} \leq r) \leq P(|X_{t^*}^{2r}| \leq r)$$

$$\leq P(|X_{t^*}^{2r} - EX_{t^*}^{2r}| \geq (1-(2a)^{-1})|EX_{t^*}^{2r}|)$$

$$\leq \frac{\text{Var } X_{t^*}^{2r}}{(1-(2a)^{-1})^2|EX_{t^*}^{2r}|^2} = \frac{K_{t^*}(2r)}{(1-(2a)^{-1})^2 M_{t^*}(2r)^2}$$

$$\leq \frac{K_t(2r)}{(1-(2a)^{-1})^2 M_t^*(2r)^2}$$

$$\leq \frac{2(1+c)M_t^*(2r)}{(1-(2a)^{-1})^2 M_t^*(2r)^2} = 2(1+c)(1-(2a)^{-1})^{-2}M_t^*(2r)^{-1}$$

$$\leq 6(1+c)^2(1-(2a)^{-1})^{-2}y_t(2r)^{-1}.$$

Here we have used Chebyshev's inequality and inequality $D_r \leq P(X_t^{2r*} \leq r)$ which implies $P(X_{t^*}^* \leq r) \leq P(X_{t^*}^{2r*} \leq r)$. Next we minimize $3(1+c^{-1})$, $3(1+c)a, 6(1+c)^2(1-(2a)^{-1})^{-2}$. Just set $3(1+c^{-1}) = 3(1+c)a = 6(1+c)^2(1-(2a)^{-1})^{-2}$. We find $a = 7/2, c = 2/7$ and $3(1+c^{-1}) = 13.5$. Thus, $D_r \leq 13.5 y_t(2r)^{-1}$. Of course, $D_r \leq \sqrt{13.5} y_t(2r)^{-1/2}$ since $D_r \leq 1$. That also covers the first case since $2\sqrt{\pi} < \sqrt{13.5}$. Applying (2.5) to $y_t(4r)^{-1/2}$ yields $D_{2r} \leq 18\sqrt{2} y_t(r)^{-1/2}$.

$Y_t = X_{t_1+t} - X_{t_1}$ with $t_1 \in [0,\infty)$ fixed is also an additive process for which the function in (2.3) equals $G_{t+t_1}(r) + K_{t+t_1}(r) - (G_{t_1}(r) + K_{t_1}(r)) + M^*_{t_1,t+t_1}(r)$, where

$$M^*_{t_1,t}(r) = \max_{t_1 \leq s \leq t} M_{t_1,s}(r),$$

$$M_{t_1,t}(r) = r^{-1}\left|B_t - B_{t_1} + \int_{1<|x|\leq r\vee 1} x\nu_{t_1,t}(dx) - \int_{r\wedge 1<|x|\leq 1} x\nu_{t_1,t}(dx)\right|,$$

for $t_1 \leq t$ and $\nu_{t_1,t} = \nu_t - \nu_{t_1}$. Since $M^*_{t_1,t_2} \geq M^*_{t_2}(r) - M^*_{t_1}(r)$ for $t_1 \leq t_2$, $P(\sup_{t_1 \leq s \leq t_2}|X_s - X_{t_1}| \leq 2r) \leq 18\sqrt{2}(y_{t_2}(r) - y_{t_1}(r))^{-1/2}$ for $t_1 < t_2$ by the result that $D_{2r} \leq 18\sqrt{2} y_t(r)^{-1/2}$. Since $y_t(r)$ is nondecreasing continuous in $t$, there are points $0 < t_1 < t_2 < \cdots < t_{k-1} < t$ such that $y_t(r)/k = y_{t_2}(r) - y_{t_1}(r) = \cdots = y_t(r) - y_{t_{k-1}}(r)$. By independence, $P(X^*_t \leq r) \leq P(X^*_{t_1} \leq r)P(\sup_{t_1\leq s\leq t_2}|X_s - X_{t_1}| \leq 2r)\cdots P(\sup_{t_{k-1}\leq s\leq t}|X_s - X_{t_{k-1}}| \leq 2r) \leq (18\sqrt{2})^k(y_{t_1}(r)(y_{t_2}(r) - y_{t_1}(r))\cdots(y_t(r) - y_{t_{k-1}}(r)))^{-1/2} = (18\sqrt{2k})^k y_t(r)^{-k/2}$. Equation (2.6) in $d=1$ has been proved.

For $d > 1$, we have $P(X^*_t \geq r) \leq \sum_{j=1}^d P(X^{(j)*}_t \geq r/d) \leq a\sum_{j=1}^d y^{(j)}_t(r/d) \leq a(3d^2)\sum_{j=1}^d y^{(j)}_t(r) = \pi_d y_t(r)$ by (2.5) and $P(X^*_t \leq r) \leq P(\max_{1\leq j\leq d}\{X^{(j)*}_t\} \leq r) \leq \min_{1\leq j\leq d}\{P(X^{(j)*}_t \leq r)\} \leq \min_{1\leq j\leq d}\{(18\sqrt{2k})^k y^{(j)}_t(r)^{-k/2}\} \leq (18\sqrt{2dk})^k(\sum_{j=1}^d y^{(j)}_t(r))^{-k/2} = A_k(d)y_t(r)^{-k/2}$. □

Let $X_t$ be an additive process in $\mathbb{R}^d$. There exists $\bar{t} \in [0,\infty]$ such that $\int_{|x|\leq 1}|x|\nu_t(dx) < \infty$ for $t \in [0,\bar{t}]$ and $\int_{|x|\leq 1}|x|\nu_t(dx) = \infty$ for $t > \bar{t}$. (E.g., $\nu_t = f(t)\nu_1$ for $t \in [0,\bar{t}]$ and $\nu_t = f(\bar{t})\nu_1 + (f(t) - f(\bar{t}))\nu_2$ for $t > \bar{t}$ where $\int_{|x|\leq 1}|x|\nu_1(dx) < \infty$, $\int_{|x|\leq 1}|x|\nu_2(dx) = \infty$ and $f$ is strictly increasing.) $\int_{|x|\leq 1}|x|\nu_t(dx)$ is a nondecreasing continuous function on $[0,\bar{t}]$. The continuous function $\gamma_0(t) = B_t - \int_{|x|\leq 1}x\nu_t(dx), t \in [0,\bar{t}]$, is called the drift of $X_t$. If $\gamma_0(t) = 0$, $M_t(r) = r^{-1}|\int_{|x|\leq r}x\nu_t(dx)|$. Let $\gamma^*_0(t) = \max_{0\leq s\leq t}|\gamma_0(s)|$. If $X_t$ is a process with additive components, the drift and its maximum for the $j$th component $X^{(j)}_t$ are denoted by $\gamma^{(j)}_0(t)$ and $\gamma^{(j)*}_0(t)$, respectively. If $X^{(j)}_t$ is monotone on $[0,\varepsilon]$, $|\gamma^{(j)}_0(t)|$ is nondecreasing on $[0,\varepsilon]$ and hence $\gamma^{(j)*}_0(t) = |\gamma^{(j)}_0(t)|$. $X_t$ is said to be drift-free initially if whenever $\int_{|x|\leq 1}|x|\nu^{(j)}_\varepsilon(dx) < \infty$ for some $\varepsilon > 0$, there exists $\varepsilon_1 \in (0,\varepsilon)$ such that $\gamma^{(j)*}_0(\varepsilon_1) = 0$. $p(t) = \nu_t(\mathbb{R}^d)$ is also a nondecreasing continuous function whenever it is finite, and there exists $\hat{t} \in [0,\infty]$ such that $p(t) < \infty$ for $t \in [0,\hat{t}]$ and $p(t) = \infty$ for $t > \hat{t}$. Recall that $X_t$ is a step process on $[0,\varepsilon], \varepsilon \leq \hat{t}$, so are its components, if and only if $\sum_{j=1}^d C^{(j)}_\varepsilon = 0$ and $\gamma^*_0(\varepsilon) = 0$. In that case for $t \in$

$[0, \varepsilon], r > 0$, $y_t(r) \leq \sum_{j=1}^{d}(\int_{|x| \leq r} |\frac{x}{r}| \nu_t^{(j)}(dx) + \int_{|x| \leq r}((\frac{x}{r})^2 \wedge 1)\nu_t^{(j)}(dx)) \leq 2p(t)$ where $p(t) = \sum_{j=1}^{d} p^{(j)}(t), p^{(j)}(t) = \nu_t^{(j)}(\mathbb{R})$. Define $G_t(r) = \sum_{j=1}^{d} G_t^{(j)}(r)$ and $G_t(r) = \nu_t(\{x \in \mathbb{R}^d : |x| > r\})$ if $X_t$ is additive in $\mathbb{R}^d$. Note that $\nu_t(\{x \in \mathbb{R}^d : |x| > r\}) \approx \sum_{i=1}^{d} G_t^{(i)}(r)$ where the constants in $\approx$ depend only on $d$.

LEMMA 2.2. *Let $X_t$ be any process with additive components.*

(i) *If $X_t$ is a step process on an interval $[0, \varepsilon]$, then for all $t \in [0, \varepsilon]$,* $\lim_{r \to 0} y_t(r) = p(t)$. *Otherwise,* $\lim_{r \to 0} y_t(r) = \infty$ *for all $t > 0$.*

(ii) $\lim_{r \to 0} r^2 y_t(r) = \sum_{j=1}^{d} C_t^{(j)}$.

(iii) *If for some $\varepsilon > 0$, $\sum_{j=1}^{d} C_\varepsilon^{(j)} = 0$, $\sum_{j=1}^{d} \int_{|x| \leq 1} |x| \nu_\varepsilon^{(j)}(dx) < \infty$, then for every $t \in (0, \varepsilon]$,* $\lim_{r \to 0} r y_t(r) = \sum_{j=1}^{d} \gamma_0^{(j)*}(t)$.

(iv) *If $X_t$ is drift-free initially and $\sum_{j=1}^{d} C_\varepsilon^{(j)} = 0$ for some $\varepsilon > 0$, then there exists $b > 0$ such that for every $t \in (0, b], \eta > 0, \lim_{r \to 0} r^\eta G_t(r) = 0$ implies $\lim_{r \to 0} r^\eta y_t(r) = 0$.*

Lemma 2.2 is standard. We omit the proof. There are also results for $r \uparrow \infty$ analogous to (i), (ii), (iii) of Lemma 2.2: (a) $\lim_{r \to \infty} y_t(r) = 0$ for all $t \geq 0$. Assume $d = 1$ below. (b) If $X_t \in L^2$, equivalently $\int_{|x| > 1} x^2 \nu_t(dx) < \infty$, and $EX_s = 0$, $s \in [0, t]$, that is, $X_s$ is in $L^2$ and centered up to time $t$, then $\lim_{r \to \infty} r^2 y_t(r) = EX_t^2 = C_t + \int x^2 \nu_t(dx)$. (c) If $X_t \in L^1$, equivalently $\int_{|x| > 1} |x| \nu_t(dx) < \infty$, then $\lim_{r \to \infty} r y_t(r) = \max_{0 \leq s \leq t} |B_s + \int_{|x| > 1} x \nu_s(dx)| = \max_{0 \leq s \leq t} |EX_s|$.

If $X_t$ is increasing on $[0, \varepsilon]$, then $C_t = 0$, $\nu_t$ has no mass on $(-\infty, 0]$ with $\int_{x \leq 1} x \nu_t(dx) < \infty$, and $B_t - \int_{x \leq 1} x \nu_t(dx)$ is nondecreasing in $t \in [0, \varepsilon]$. Thus, $M_t(r) = r^{-1}(B_t - \int_{x \leq 1} x \nu_t(dx) + \int_{x \leq r} x \nu_t(dx))$ is nondecreasing in $t$ and $G_t(r) + M_t(r) = r^{-1}(B_t - \int_{x \leq 1} x \nu_t(dx)) + \int (x/r) \wedge 1 \nu_t(dx)$. It follows that $y_t(r) \leq 2\theta y_t(\theta r)$ for $\theta > 1, M_t^*(r) = M_t(r), M_t(r) \geq K_t(r), G_t(r) + M_t(r) \leq y_t(r) \leq 2(G_t(r) + M_t(r))$, and $G_t(r) + M_t(r)$ is nondecreasing in $t$ and nonincreasing continuous in $r$. For the obvious reason, we use $G_t(r) + M_t(r)$ instead of $y_t(r)$. For the Laplace transform of $X_t$, we have $Ee^{-\lambda X_t} = e^{-\psi(t,\lambda)}, \lambda > 0$, where $\psi(t, \lambda) = \lambda \gamma_0(t) + g_t(\lambda), g_t(\lambda) = \int_0^\infty (1 - e^{-\lambda x}) \nu_t(dx)$. Clearly $G_t(r) + M_t(r) = r^{-1} \gamma_0(t) + \int (x/r) \wedge 1 \nu_t(dx)$. Since $e^{-1}(x \wedge 1) < 1 - e^{-x} < x \wedge 1, x > 0$, $e^{-1} \int_0^\infty (x/r) \wedge 1 \nu_t(dx) \leq g_t(r^{-1}) \leq \int_0^\infty (x/r) \wedge 1 \nu_t(dx)$. Therefore, $y_t(r) \approx G_t(r) + M_t(r) \approx r^{-1} \gamma_0(t) + g_t(r^{-1})$. The same can be said for a decreasing process as well as any process with monotone components. If a real $X_t$ is symmetric on $[0, \varepsilon]$, that is, $E \exp\{i\lambda X_t\}$ is real, then $B_t = 0$ and $\nu_t$ is symmetric for $t \in [0, \varepsilon]$, in which case $M_t(r)$ vanishes, $y_t(r) = G_t(r) + K_t(r)$ and $y_t(r) \leq \theta^2 y_t(\theta r)$ for $\theta > 1$.

$y_t(r)$ is comparable to a function that is jointly continuous and strictly decreasing in $r$. Let

$$I_t(r) = r^{-1} \int_0^r y_t(x)^{-1} \, dx, \qquad \dot{y}_t(r) = I_t(r)^{-1}, \qquad t > 0, r > 0.$$

By (2.5), $y_t(x) \geq 2^{-1} y_t(r)$ for $x \in (0, r]$ and $I_t(r) \geq r^{-1} \int_{r/2}^r y_t(x)^{-1} \, dx \geq r^{-1} \int_{r/2}^r y_t(r/2)^{-1} \, dx = 4^{-1} y_t(r/2)^{-1} \geq 48^{-1} y_t(r)^{-1}$, which shows that for $t > 0$, $r > 0$,

$$48^{-1} \leq y_t(r)/\dot{y}_t(r) \leq 2$$

and by (2.5), for $\theta > 1, t > 0, r > 0$, $k_1 \theta^{-2} \dot{y}_t(r) \leq \dot{y}_t(\theta r) \leq k_2 \dot{y}_t(r)$, where $k_1 = 288^{-1}, k_2 = 192$. By (2.5), $2^{-1} y_t(r) \leq \inf_{0 < x \leq r} y_t(x) \leq y_t(r)$. If we use $\inf_{0 < x \leq r} y_t(x)$ instead of $y_t(r)$, $\dot{y}_t(r)$ is strictly decreasing in $r$.

LEMMA 2.3. *$\dot{y}_t(r)$ is jointly continuous.*

PROOF. $I_t(r)$ is well defined since $y_t(r)$ is rcll in $r$, nonincreasing in $t$ since $y_t(r)$ is nondecreasing in $t$, continuous in $t$ by the dominated convergence theorem since $y_t(x) \geq 2^{-1} y_t(r)$ for $x \in (0, r]$ by (2.5), and absolutely continuous in $r$ because of the way it is defined. [Hence, $\dot{y}_t(r)$ is nondecreasing continuous in $t$ and absolutely continuous in $r$.] It is enough to show that $I_t(r)$ is jointly continuous in $d = 1$. First we claim that given $r' > 0, t' > 0, \varepsilon > 0$, there exists $\delta > 0$ not depending on $r, t_1, t_2$ such that $y_{t_2}(r) - y_{t_1}(r) < \varepsilon$ whenever $r \geq r'$, $t_2 - t_1 < \delta$, $t_1, t_2 \in [0, t']$. The definition of $\nu_t$ and an approximation argument show that for any $A \in \mathcal{B}(\mathbb{R}^d)$ and Borel function $f$ satisfying $\int_A |f(x)| \nu_t(dx) < \infty$, $\int_A f(x) \nu_t(dx) = \int_{[0,t] \times A} f(x) \nu(ds, dx)$. Let $Q_t(r) = G_t(r) + K_t(r)$. Then $Q_{t_2}(r) - Q_{t_1}(r) = r^{-2}(C_{t_2} - C_{t_1}) + \int_{[t_1, t_2] \times \mathbb{R}} (x/r)^2 \wedge 1 \nu(ds, dx) \leq r'^{-2}(C_{t_2} - C_{t_1}) + \int_{[t_1, t_2] \times \mathbb{R}} (x/r')^2 \wedge 1 \nu(ds, dx) = Q_{t_2}(r') - Q_{t_1}(r') < \varepsilon$ since $Q_t(r')$ is uniformly continuous on $[0, t']$. It remains to show that $M_{t_2}^*(r) - M_{t_1}^*(r) \leq M_{t_1, t_2}^*(r) < \varepsilon$. Since $B_t$ is uniformly continuous on $[0, t'], r^{-1}|B_{t_2} - B_{t_1}| \leq r'^{-1}|B_{t_2} - B_{t_1}| < \varepsilon$. For $r < 1$, $r^{-1} \int_{r < |x| \leq 1} |x| \nu_{t_1, t_2}(dx) \leq r'^{-1} \int_{r' < |x| \leq 1} |x| \nu_{t_1, t_2}(dx) \leq r'^{-1} \int_{|x| > r'} \nu_{t_1, t_2}(dx) = r'^{-1}(\nu_{t_2}(\{x : |x| > r'\}) - \nu_{t_1}(\{x : |x| > r'\})) < \varepsilon$ since $\nu_t(\{x : |x| > r'\})$ is uniformly continuous on $[0, t']$. Similarly, for $r > 1$, $r^{-1} \int_{1 < |x| \leq r} |x| \nu_{t_1, t_2}(dx) \leq \int_{|x| > 1} \nu_{t_1, t_2}(dx) < \varepsilon$. The claim is proved. For $0 < t'' < t', t_1 < t_2, t_1, t_2 \in [t'', t'], y_{t_1}(r)^{-1} - y_{t_2}(r)^{-1} = (y_{t_1}(r) \times y_{t_2}(r))^{-1}(y_{t_2}(r) - y_{t_1}(r)) \leq (4/y_{t''}(r')^2)(y_{t_2}(r) - y_{t_1}(r))$ by (2.5). It follows from the claim above that given $r' > 0, 0 < t'' < t', \varepsilon > 0$, there exists $\delta > 0$ not depending on $r, t_1, t_2$ such that $y_{t_1}(r)^{-1} - y_{t_2}(r)^{-1} < \varepsilon$ whenever $r \geq r', t_2 - t_1 < \delta, t_1, t_2 \in [t'', t']$. Next fix a point $(t_0, r_0) \in (0, \infty) \times (0, \infty)$. Since $I_{t_0}(r)$ is (absolutely) continuous in $r$, there is $\delta_1 > 0$ such that $|I_{t_0}(r) - I_{t_0}(r_0)| < \varepsilon$ for $r \in (r_0 - \delta_1, r_0 + \delta_1)$ with $r_1 = r_0 - \delta_1 > 0$. On the other hand, by the fact that $I_t(r_1)$ is continuous in $t$ and by the result following the claim, there exists $\delta_2 > 0$ such that when $t - t_0 < \delta_2, t_0 < t, I_{t_0}(r_1) - I_t(r_1) < \varepsilon, y_{t_0}(s)^{-1} - y_t(s)^{-1} < \varepsilon$ for all $s \geq r_1$. Thus, for $r \in (r_0 - \delta_1, r_0 + \delta_1), t - t_0 < \delta_2, t_0 < t, |I_t(r) - I_{t_0}(r_0)| \leq |I_{t_0}(r) - I_{t_0}(r_0)| + |I_t(r) - I_{t_0}(r)| < \varepsilon + I_{t_0}(r) - I_t(r)$ and

$$
\begin{aligned}
I_{t_0}(r) &- I_t(r) \\
&= r^{-1} \int_0^r (y_{t_0}(s)^{-1} - y_t(s)^{-1}) \, ds \\
&\leq r_1^{-1} \int_0^{r_1} (y_{t_0}(s)^{-1} - y_t(s)^{-1}) \, ds + r^{-1} \int_{r_1}^r (y_{t_0}(s)^{-1} - y_t(s)^{-1}) \, ds \\
&= I_{t_0}(r_1) - I_t(r_1) + r^{-1} \int_{r_1}^r (y_{t_0}(s)^{-1} - y_t(s)^{-1}) \, ds \\
&< \varepsilon + r^{-1}(r - r_1)\varepsilon < 2\varepsilon.
\end{aligned}
$$

The treatment for $t_0 - t < \delta_2, t < t_0$ is completely analogous. $\square$

**3. The quasiconvex function method.** A sequence $\sigma_n \downarrow 0$ is called the $\Sigma$-sequence if $\frac{\sigma_{n-1}}{\sigma_n} \cdot \sigma_n^\eta \to 0$ as $n \to \infty$ for all $\eta > 0$, which implies that $(\frac{\sigma_{n-1}}{\sigma_n})^\varepsilon \cdot \sigma_n^\eta \to 0$ and $(\frac{\sigma_{n-1}}{\sigma_n})^\varepsilon \cdot \sigma_{n-1}^\eta \to 0$ for all $\varepsilon > 0$. Some of the $\Sigma$-sequences are constructed from continuous slow functions $c$. If $\sigma_{n+1}/\sigma_n \geq c(\sigma_{n+1})$, $\sigma_n \downarrow 0$ is a $\Sigma$-sequence since $(\sigma_{n-1}/\sigma_n) \cdot \sigma_n^\eta \leq \sigma_n^\eta/c(\sigma_n) \to 0$. [For any $s_n \in (0,1)$, there is $s_{n+1} < s_n$ such that $s_{n+1}/c(s_{n+1}) = s_n$ because $s_n/c(s_n) > s_n$ and $t/c(t) \to 0$. $\lim_{n\to\infty} s_n = 0$ holds also.] Let $X_t$ be a process in $\mathbb{R}^d$ continuous in probability with $X_0 = 0$ and $v$ a nondecreasing function. Define

$$\underline{\delta} = \inf\{\eta > 0 : P(X^*_{v(t_n)} \leq t_n^{1/\eta} \text{ i.o.}) = 1 \text{ for some sequence } t_n \downarrow 0\},$$

$$\overline{\delta} = \sup\{\eta > 0 : P(X^*_{v(\sigma_n)} \leq \sigma_n^{1/\eta} \text{ i.o.}) = 0 \text{ for some } \Sigma\text{-sequence } \sigma_n \downarrow 0\}.$$

(Both the sequence $t_n \downarrow 0$ and the $\Sigma$-sequence $\sigma_n \downarrow 0$ in braces depend on $\eta$.) If $\eta > \underline{\delta}$, the stronger result that $\liminf_{n\to\infty} t_n^{-1/\eta} X^*_{v(t_n)} = 0$ a.s. for some sequence $t_n \downarrow 0$ holds (which implies that $\liminf_{t\to 0} t^{-1/\eta} X^*_{v(t)} = 0$ a.s.). If $\eta < \overline{\delta}$, there exists $\eta_1 > \eta$ such that $X^*_{v(\sigma_n)} > \sigma_n^{1/\eta_1}$ for all large $n$. Therefore, for $t \in [\sigma_{n+1}, \sigma_n]$, $X^*_{v(t)}/t^{1/\eta} \geq X^*_{v(\sigma_{n+1})}/\sigma_n^{1/\eta} > \sigma_{n+1}^{1/\eta_1}/\sigma_n^{1/\eta} = (\sigma_{n+1}/\sigma_n)^{1/\eta_1}/\sigma_n^{(1/\eta - 1/\eta_1)} \to \infty$, which implies that $\lim_{t\to 0} t^{-1/\eta} X^*_{v(t)} = \infty$ a.s. If we define

$$\dot{\delta}_1 = \inf\left\{\eta > 0 : \liminf_{r\to 0} P(X^*_{v(r)} \geq r^{1/\eta}) = 0\right\},$$

$$\dot{\delta}_2 = \sup\left\{\eta > 0 : \sum P(X^*_{v(\sigma_n)} \leq \sigma_n^{1/\eta}) < \infty \text{ for some } \Sigma\text{-sequence } \sigma_n \downarrow 0\right\},$$

the Borel–Cantelli lemma and Fatou's lemma imply that $\dot{\delta}_2 \leq \overline{\delta}$ and $\underline{\delta} \leq \dot{\delta}_1$. ($\overline{\delta} \leq \underline{\delta}$.) By the same token, if we define

$$\underline{\beta} = \sup\{\eta > 0 : P(X^*_{v(t_n)} \geq t_n^{1/\eta} \text{ i.o.}) = 1 \text{ for some sequence } t_n \downarrow 0\},$$

$$\overline{\beta} = \inf\{\eta > 0 : P(X^*_{v(\sigma_n)} \geq \sigma_n^{1/\eta} \text{ i.o.}) = 0 \text{ for some } \Sigma\text{-sequence } \sigma_n \downarrow 0\},$$

$$\dot{\beta}_1 = \sup\left\{\eta > 0 : \liminf_{r \to 0} P(X^*_{v(r)} \leq r^{1/\eta}) = 0\right\},$$

$$\dot{\beta}_2 = \inf\left\{\eta > 0 : \sum P(X^*_{v(\sigma_n)} \geq \sigma_n^{1/\eta}) < \infty \text{ for some } \Sigma\text{-sequence } \sigma_n \downarrow 0\right\},$$

we have $\limsup_{t \to 0} t^{-1/\eta} X^*_{v(t)} = \infty$ a.s. for $\eta < \underline{\beta}$ and $\lim_{t \to 0} t^{-1/\eta} X^*_{v(t)} = 0$ a.s. for $\eta > \overline{\beta}$ while $\dot{\beta}_1 \leq \underline{\beta} \leq \overline{\beta} \leq \dot{\beta}_2$. Clearly, $\dot{\delta}_2 \leq \dot{\beta}_1, \dot{\delta}_1 \leq \dot{\beta}_2$. To define $v$ with $\dot{\beta}_2 < \infty$, we fix a number $\kappa \in (0, \infty)$ along with a $\Sigma$-sequence $\bar{\sigma}_n \downarrow 0$. There is always a sequence $v_n \downarrow 0$ such that $\sum P(X^*_{v_n} \geq \bar{\sigma}_n^{1/\kappa}) < \infty$. Let $v$ be a nondecreasing function with values $v_n$ at $\bar{\sigma}_n$. Then $\dot{\beta}_2 \leq \kappa$.

We cannot get anything better than $\overline{\delta}, \underline{\delta}, \underline{\beta}, \overline{\beta}$ if $v$ is fixed. If $\liminf_{n \to \infty} t_n^{-1/\eta} \times X^*_{v(t_n)} = 0$ a.s. for some sequence $t_n \downarrow 0$, $\eta \geq \underline{\delta}$. If $\liminf_{t \to 0} t^{-1/\eta} X^*_{v(t)} = 0$ a.s., then with respect to each $\omega \in \Omega$, there is a sequence $t_n^\omega \downarrow 0$ such that $X^{\omega*}_{v(t_n^\omega)}/(t_n^\omega)^{1/\eta} \leq 1$ for large $n$ except for $\omega$ in a $P$-null set. The sequences can be extracted technically from a fixed deterministic sequence $t_n \downarrow 0$, that is, $P(X^*_{v(t_n)} \leq t_n^{1/\eta} \text{ i.o.}) = 1$. Thus, $\eta \geq \underline{\delta}$. If $\lim_{t \to 0} t^{-1/\eta} X^*_{v(t)} = \infty$ a.s., then for all sequences $s_n \downarrow 0$, $X^*_{v(s_n)}/s_n^{1/\eta} \geq 1$ for large $n$ a.s., that is, $P(X^*_{v(s_n)} \leq s_n^{1/\eta} \text{i.o.}) = 0$. Hence, $\eta \leq \overline{\delta}$. Same goes for $\overline{\beta}, \underline{\beta}$.

Let $X_t$ be any process with additive components. Define

$$\delta_1 = \inf\left\{\eta > 0 : \liminf_{r \to 0} y_{v(r)}(r^{1/\eta}) = 0\right\},$$

$$\delta_2 = \sup\left\{\eta > 0 : \sum y_{v(\sigma_n)}(\sigma_n^{1/\eta})^{-1} < \infty \text{ for some } \Sigma\text{-sequence } \sigma_n \downarrow 0\right\},$$

$$\beta_1 = \sup\left\{\eta > 0 : \liminf_{r \to 0} y_{v(r)}(r^{1/\eta})^{-1} = 0\right\},$$

$$\beta_2 = \inf\left\{\eta > 0 : \sum y_{v(\sigma_n)}(\sigma_n^{1/\eta}) < \infty \text{ for some } \Sigma\text{-sequence } \sigma_n \downarrow 0\right\}.$$

Clearly, $\delta_2 \leq \beta_1, \delta_1 \leq \beta_2$. By Lemma 2.1, $\delta_2 \leq \dot{\delta}_2$, $\dot{\delta}_1 \leq \delta_1$, $\beta_1 \leq \dot{\beta}_1, \dot{\beta}_2 \leq \beta_2$. Similarly, to define $v$ with $\beta_2 < \infty$, we can preselect a number $\kappa \in (0, \infty)$ and a $\Sigma$-sequence $\bar{\sigma}_n \downarrow 0$. Since $y_\varepsilon(r) \downarrow 0$ as $\varepsilon \downarrow 0$, there is a sequence $v_n \downarrow 0$ such that $\sum y_{v_n}(\bar{\sigma}_n^{1/\kappa}) < \infty$. If $v$ is a nondecreasing function taking values $v_n$ at $\bar{\sigma}_n$, then $\beta_2 \leq \kappa$. Of course, if $\beta_2 < \infty$, $v$ has to be defined in this way. We wish the definition of $v$ given above to be more specific. $v$ should have the information about the case $\delta_1 = \delta_2, \beta_1 = \beta_2$ and should be able to equal $t$ when $X_t$ is a Lévy process. We define $v$ as follows. Select a quasiconvex function $\phi$ and a constant $b \in (0, \infty)$. Equation $y_{v(t)}(b) = \phi(t)$ defines a nondecreasing function $v$ with continuous inverse $u$ since $y_t(r)$ is nondecreasing continuous in $t$ with $y_0(r) = 0$. For example, $u(t) = y_t(b)^{1/p}$ for $\phi(t) = t^p, p > 0$ while $v(t) = u(t) = t$ for $\phi(t) = h(b)t$ in the case of Lévy processes. By (2.5) and quasiconvexity of $\phi$, $y_{v(r)}(r^{1/\eta}) \le c y_{v(r)}(b) r^{-2/\eta} = c\phi(r) r^{-2/\eta} \le c_1 r^\sigma r^{-2/\eta}, r < 1 \wedge b$. Taking a $\Sigma$-sequence such as $\sigma_n = 2^{-n}$ shows that $\beta_2 \le 2/\sigma$. The result in (1.1) remains unchanged when $u(t)^{-1/\eta} X_t^*$ is replaced by $t^{-1/\eta} X_{v(t)}^*$. We have proved

THEOREM 3.1. *Let $X_t$ be any process in $\mathbb{R}^d$ with additive components and $\delta_1, \delta_2, \beta_1, \beta_2, u$ as given above. Then (1.1) holds.*

$y_{v(t)}(b) = \phi(t)$ is quasiconvex. In fact, for all $r \in (0, b)$, $y_{v(t)}(r)$ is quasiconvex as well with $y_{v(t_2)}(r)/y_{v(t_1)}(r) \ge \rho_r (t_2/t_1)^\sigma, t_1 < t_2$, where $\rho_r \ge 6^{-1}(r/b)^2 \rho$ by (2.5). But $6^{-1}(r/b)^2 \rho$ is not a slow function. $y_{v(t)}(r)$ is called quasiconvex with respect to a nondecreasing function $v > 0$ if

$$(3.1) \quad y_{v(t_2)}(r)/y_{v(t_1)}(r) \ge c(r)(t_2/t_1)^\sigma, \qquad 0 < t_1 < t_2 \le t_0, r \in (0, b],$$

with a slow function $c(r)$ and a constant $\sigma > 0$. Equation (3.1) means $\rho_r \ge c(r)$. Conversely, any $v$ satisfying (3.1) is valid for Theorem 3.1 since $y_{v(t)}(b)$ is quasiconvex. Define

$$n(r) = \inf\{t > 0 : y_{v(t)}(r) > m\}$$

with $m \in (0, 2^{-1} y_{v(t_0)}(b) \wedge (aK(d))^{-1})$ held fixed. Equation (2.5) implies that $n(r)$ is a finitely determined function. Let

$$\delta = \inf\left\{\eta \ge 0 : \liminf_{r \to 0} r^\eta n(r)^{-1} = 0\right\},$$

$$\beta = \inf\left\{\eta \ge 0 : \lim_{r \to 0} r^\eta n(r)^{-1} = 0\right\}$$

and for any fixed constant $l \in (0, t_0]$,

$$\delta_P = \sup\left\{\eta \ge 0 : \limsup_{r \to 0} r^{-\eta} \int_0^l P(X_{v(t)}^* \le r)\, dt < \infty\right\},$$

$$\beta_P = \sup\left\{\eta \ge 0 : \liminf_{r \to 0} r^{-\eta} \int_0^l P(X_{v(t)}^* \le r)\, dt < \infty\right\},$$

$$\delta_E = \sup\left\{\eta \ge 0 : \int_0^l E(X_{v(t)}^*)^{-\eta}\, dt < \infty\right\}.$$

Note that $\int_0^l P(X_{v(t)}^* \le r)\, dt = E(T_r^v \wedge l)$ where $T_r^v = \inf\{t > 0 : |X_{v(t)}| > r\}$. In the case of Lévy processes ($v(t) = t$), $n(r) = mh(r)^{-1}$. $n(r)$ is an $h(r)^{-1}$ analogy.

THEOREM 3.2. *In Theorem 3.1 if $y_{v(t)}(r)$ is quasiconvex, then $\delta_1 = \delta_2 = \delta = \delta_P = \delta_E, \beta_1 = \beta_2 = \beta = \beta_P$. $\beta \le 1/\sigma$ if $\sum_{j=1}^d \int_{|x|\le 1} |x| \nu_\tau^{(j)}(dx) < \infty$ and $\sum_{j=1}^d C_\tau^{(j)} = 0$ for some $\tau > 0$. If $X_t$ is drift-free initially, $\sum_{j=1}^d C_\tau^{(j)} = 0$ for some $\tau > 0$, and $G_{v(t_2)}(r)/G_{v(t_1)}(r) \le c'(r)^{-1}(t_2/t_1)^{\sigma'}$ for $0 < t_1 < t_2 \le t_0'$, $r \in (0, r_0)$ with a constant $\sigma' \le \sigma$ and a slow function $c'$, then $\beta = \inf\{\eta \ge 0 : \lim_{r \to 0} r^\eta \bar{n}(r)^{-1} = 0\}$ where $\bar{n}(r) = \inf\{t > 0 : G_{v(t)}(r) > m\}$. If $X_t$ is drift-free initially with increasing components, then $\delta = \sup\{\eta \ge 0 : \lim_{r \to \infty} r^\eta \hat{n}(r) = 0\}$, $\beta = \inf\{\eta \ge 0 : \lim_{r \to \infty} r^\eta \hat{n}(r) = \infty\}$, where $\hat{n}(r) = \inf\{t > 0 : g_{v(t)}(r) > m\}$ with $g_t(r) = \sum_{i=1}^d g_t^{(i)}(r)$.*

PROOF. (i) $\delta_1 = \delta_2, \beta_1 = \beta_2$: Define $g_\varepsilon(r) = r^{\sigma\varepsilon/2}$, $\delta_\varepsilon = \inf\{\eta > 0 : \liminf_{r \to 0} g_\varepsilon(r) y_{v(r)}(r^{1/\eta}) = 0\}$, and $\beta_\varepsilon = \sup\{\eta > 0 : \liminf_{r \to 0} g_\varepsilon(r) \times y_{v(r)}(r^{1/\eta})^{-1} = 0\}$. Notice that $\delta_\varepsilon \uparrow \delta^* \le \delta_2, \beta_\varepsilon \downarrow \beta^* \ge \beta_2$, as $\varepsilon \downarrow 0$. If $\delta^* < \delta_1$, pick $\eta \in (\delta^*, \delta_1)$ and $\varepsilon \in (0, \delta_1/\eta - 1)$. Since $\delta_\varepsilon < \eta$, there exists an $\eta_1 < \eta$ such that $g_\varepsilon(r_n) y_{v(r_n)}(r_n^{1/\eta_1}) \to 0$ for some sequence $r_n \downarrow 0$. Let $t_n = r_n^{1+\varepsilon}$. Then $t_n^{1/(1+\varepsilon)\eta_1} = r_n^{1/\eta_1}$. Since $y_{v(t)}(r)$ is quasiconvex and since $c(r^{1/\eta})^{-1} \le c_\eta r^{-\sigma\varepsilon/2}$, $y_{v(t_n)}(t_n^{1/(1+\varepsilon)\eta_1})/y_{v(r_n)}(r_n^{1/\eta_1}) = y_{v(r_n^{1+\varepsilon})}(r_n^{1/\eta_1})/y_{v(r_n)}(r_n^{1/\eta_1}) \le c(r_n^{1/\eta_1})^{-1}(r_n^{1+\varepsilon}/r_n)^\sigma \le c_{\eta_1} r_n^{-\sigma\varepsilon/2} r_n^{\sigma\varepsilon} = c_{\eta_1} g_\varepsilon(r_n)$. Since $(1+\varepsilon)\eta_1 < \delta_1$, $y_{v(t_n)}(t_n^{1/(1+\varepsilon)\eta_1}) \ge c > 0$. Thus, $g_\varepsilon(r_n) y_{v(r_n)}(r_n^{1/\eta_1}) \ge c/c_{\eta_1} > 0$. We have a contradiction. The argument for $\beta^* = \beta_1$ is similar.

(ii) $\delta_1 = \delta, \beta_1 = \beta$: If $\delta_1 < \delta, r^\eta n(r)^{-1} \ge c > 0$ for any $\eta \in (\delta_1, \delta)$ and hence $r^{\eta_1} \ge n(r)$ for $\eta_1 \in (\delta_1, \eta)$ and $r$ small. Since $\delta_1 < \eta_1$, there exists an $\eta_2 \in (\delta_1, \eta_1)$ such that $\liminf_{r \to 0} y_{v(r)}(r^{1/\eta_2}) = 0$. Since $r^{\eta_2} \ge n(r)$ as well, by quasiconvexity and the facts that $c(r) \ge r^\varepsilon$ for all $\varepsilon > 0$ and $y_{v(n(r))}(r) = m$, $y_{v(r^{\eta_2})}(r)/m = y_{v(r^{\eta_2})}(r)/y_{v(n(r))}(r) \ge c(r)(r^{\eta_2}/n(r))^\sigma$, which implies that $\liminf_{r \to 0} r^{\eta_1} n(r)^{-1} = 0$ contradicting $\eta_1 < \delta$. If $\delta < \delta_1, r_n^\eta n(r_n)^{-1} \to 0$ for some $\eta \in (\delta, \delta_1)$ and a sequence $r_n \downarrow 0$, which implies that $r_n^\eta \le n(r_n)$ and $r_n^{\eta_1} \le n(r_n)$ for $\eta_1 \in (\eta, \delta_1)$. By quasiconvexity, $m/y_{v(r_n^{\eta_1})}(r_n) \ge c(r_n)(n(r_n)/r_n^{\eta_1})^\sigma \ge (n(r_n)/r_n^\eta)^\sigma$. Thus, $y_{v(r_n^{\eta_1})}(r_n) \to 0$ and $\delta_1 \le \eta_1$. That is a contradiction. If $\beta_1 < \beta$, then for any $\eta \in (\beta_1, \beta)$ there exists a sequence $r_n \downarrow 0$ such that $r_n^{-\eta} n(r_n) \to 0$. So, $r_n^\eta \ge n(r_n)$ and $r_n^{\eta_1} \ge n(r_n)$ for $\eta_1 \in (\beta_1, \eta)$. By quasiconvexity, $y_{v(r_n^{\eta_1})}(r_n)/m \ge c(r_n)(r_n^{\eta_1}/n(r_n))^\sigma \ge (r_n^\eta/n(r_n))^\sigma$. Thus, $y_{v(r_n^{\eta_1})}(r_n)^{-1} \to 0$ and $\eta_1 \le \beta_1$ contradicting $\eta_1 > \beta_1$. Lastly, if $\beta < \beta_1$, then for any $\eta \in (\beta, \beta_1)$, $r^\eta n(r)^{-1} \to 0$; that is, $n(r) \ge r^\eta$. Since $\eta < \beta_1$, there is an $\eta_1 \in (\eta, \beta_1)$ and a sequence $r_n \downarrow 0$ such that $y_{v(r_n^{\eta_1})}(r_n)^{-1} \to 0$. But $n(r) \ge r^{\eta_1}$, so by quasiconvexity we have $y_{v(r^{\eta_1})}(r) \le m(r^\eta/n(r))^\sigma \to 0$. That is a contradiction.

(iii) $\delta = \delta_P, \beta = \beta_P$ and $\beta \le 1/\sigma$ with the condition as stated: First we prove that if $X_t$ is a step process initially, then $\inf_{r > 0} n(r) > 0$ and otherwise $\lim_{r \to 0} n(r) = 0$. Let $X_t$ be a step process up to time $\hat{t}$. Fix $r > 0$. If $y_\infty(r) \leq m$, then $n(r) = \infty$. Suppose that $y_{v(t^*)}(r) > m$ for some $t^* > 0$. Then $n(r) < \infty$. If $v(n(r)) < \hat{t}$, $p(v(n(r))) \geq 2^{-1} y_{v(n(r))}(r) = 2^{-1}m$ where $p(t) = \sum_{j=1}^{d} p^{(j)}(t), p^{(j)}(t) = \nu_t^{(j)}(\mathbb{R})$. Thus, there exists a positive constant $K$ such that $n(r) \geq K$ for all $r > 0$. In the second case, since $\lim_{r \to 0} y_{v(t)}(r) = \infty$ for any $t \in (0, \infty)$ by Lemma 2.2(i), there exists $r_0 > 0$ depending on $t$ such that $y_{v(t)}(r_0) > 2m$. Therefore, for $r < r_0, y_{v(t)}(r) \geq 2^{-1} y_{v(t)}(r_0) > m$ and $n(r) \leq t$. If $\lim_{r \to 0} n(r) = 0$ fails, there is a sequence $r_n \downarrow 0, r_n \in (0, r_0)$ such that $n(r_n) \geq \delta$ for some $\delta > 0$ and hence $y_{v(\delta)}(r_n) \leq y_{v(n(r_n))}(r_n) = m < \infty$ contradicting $\lim_{r \to 0} y_{v(\delta)}(r) = \infty$ according to Lemma 2.2(i). Let us move on to prove (iii) in this case. Choose $r_0$ as above such that $n(r) \leq l$ for all $r \in (0, r_0)$. Note that for $t \in (0, n(r)], y_{v(t)}(r) \leq m$ and hence by the first bound in (2.6), $P(X_{v(t)}^* \geq r) \leq aK(d) y_{v(t)}(r) \leq aK(d)m$. It follows that $\int_0^l P(X_{v(t)}^* \leq r) \, dt \geq \int_0^{n(r)} P(X_{v(t)}^* \leq r) \, dt = \int_0^{n(r)} (1 - P(X_{v(t)}^* > r)) \, dt \geq (1 - aK(d)m)n(r)$. For $t \in [n(r), l]$, $y_{v(t)}(r)/m \geq c(r)(t/n(r))^\sigma$ by quasiconvexity. Choose an integer $k > 2/\sigma$ and let $\theta = k\sigma/2 - 1 > 0$. By the second bound in (2.6), $P(X_{v(t)}^* \leq r) \leq c(y_{v(t)}(r))^{-k/2} \leq c(m^{-1}c(r)^{-1})^{k/2}(n(r)/t)^{1+\theta} = c_1(r)^{-1}(n(r)/t)^{1+\theta}$ where $c_1(r)$ is also a slow function. Thus, $\int_0^l P(X_{v(t)}^* \leq r) \, dt = \int_0^{n(r)} + \int_{n(r)}^l \leq n(r) + \int_{n(r)}^l c_1(r)^{-1}(n(r)/t)^{1+\theta} \, dt = (1 + (\theta c_1(r))^{-1})n(r) - (l^\theta \theta c_1(r))^{-1} n(r)^{1+\theta} < c_2(r)^{-1} n(r)$ where $c_2(r)$ is another slow function. Hence $\delta = \delta_P, \beta = \beta_P$. Choose $t \leq t_0$ such that $v(t) \leq \tau$ and then choose $r_0$ such that $n(r) \leq t$ for all $r \in (0, r_0)$. By quasiconvexity, $r^{1/\sigma} n(r)^{-1} \leq t^{-1} m^{-1/\sigma} c(r)^{-1/\sigma} (r y_{v(t)}(r))^{1/\sigma}$ and by Lemma 2.2(iii), $\lim_{r \to 0} r y_{v(t)}(r) = \sum_{j=1}^d \gamma_0^{(j)*}(v(t)) < \infty$. Thus, $\beta \leq 1/\sigma$. Assume that $X_t$ is a step process initially. Then $\beta = 0$ ($< 1/\sigma$) since $\inf_{r > 0} n(r) > 0$. By choosing $m \in (0, 2^{-1} y_{v(l)}(r_0) \wedge (aK(d))^{-1})$ for any $r_0 > 0$, we redefine $n(r)$. Then $n(r) \leq l$ for all $r \in (0, r_0)$ and $\int_0^l P(X_{v(t)}^* \leq r) \, dt \geq (1 - aK(d)m)n(r)$, which shows that $\inf_{r > 0} \int_0^l P(X_{v(t)}^* \leq r) \, dt > 0$. Thus, $\delta_P \leq \beta_P = 0$.

(iv) $\delta_P = \delta_E$: Let $\zeta_t$ be a process taking nonnegative values. Define for $l > 0, r_0 > 0, g(r) = \int_0^l P(\zeta_t \leq r) \, dt, r \in (0, r_0], \delta' = \sup\{\eta \geq 0 : \int_0^l E\zeta_t^{-\eta} \, dt < \infty\}$, $\delta'' = \sup\{\eta \geq 0 : \sup_{0 < r \leq r_0} r^{-\eta} g(r) < \infty\} = \sup\{\eta \geq 0 : \limsup_{r \to 0} r^{-\eta} g(r) < \infty\}$. Then $\delta' = \delta''$. This is not difficult to prove. Clearly $g(r)$ is nondecreasing and bounded by $l$. For $\eta > 0$, $E\zeta_t^{-\eta} = \eta \int_0^\infty x^{-\eta-1} P(\zeta_t \leq x) \, dx$. Therefore, $\int_0^l E\zeta_t^{-\eta} \, dt = \eta \int_0^{r_0} x^{-\eta-1} g(x) \, dx + \eta \int_{r_0}^\infty x^{-\eta-1} (\int_0^l P(\zeta_t \leq x) \, dt) \, dx$, which shows that $\int_0^l E\zeta_t^{-\eta} \, dt < \infty$ if and only if $\int_0^{r_0} x^{-\eta-1} g(x) \, dx < \infty$. For $r \leq r_0/2, \int_0^{r_0} x^{-\eta-1} g(x) \, dx \geq \int_r^{2r} x^{-\eta-1} g(x) \, dx \geq g(r) \int_r^{2r} x^{-\eta-1} \, dx = kr^{-\eta} g(r)$, $k = \eta^{-1}(1 - 2^{-\eta})$. It follows that $\delta' \leq \delta''$. If $\eta < \delta''$, $\sup_{0 < r \leq r_0} r^{-\eta_1} g(r) < \infty$ for some $\eta_1 > \eta$. Thus $\int_0^{r_0} x^{-\eta-1} g(x) \, dx = \int_0^{r_0} x^{\eta_1 - \eta - 1} x^{-\eta_1} g(x) \, dx \leq k_1 \sup_{0 < r \leq r_0} r^{-\eta_1} g(r) < \infty$ where $k_1 = \int_0^{r_0} x^{\eta_1 - \eta - 1} \, dx < \infty$, which implies that $\delta'' \leq \delta'$. Take $\zeta_t = X_{v(t)}^*$ to finish.

(v) The last two statements in the theorem: Let $\bar{\beta} = \inf\{\eta \geq 0 : \lim_{r \to 0} r^\eta \times \bar{n}(r)^{-1} = 0\}$. Since $G_t(r) \leq y_t(r)$, $\bar{\beta} \leq \beta$. We prove the opposite. We may assume that $\sigma' = \sigma$ and yet we may also assume that $X_t$ is not a step process initially, for otherwise $\bar{\beta} = \beta = 0$. Let $b$ be the constant in Lemma 2.2(iv). Choose $t \leq t_0 \wedge t_0'$ such that $v(t) \leq b$ and then choose $r_0$ such that $n(r) \leq t$ for all $r \in (0, r_0)$. Note that $G_{v(\bar{n}(r))}(r) = m = y_{v(n(r))}(r)$. It follows that $G_{v(t)}(r)^{1/\sigma} \leq c_1(r)^{-1}/\bar{n}(r)$ by the condition on $G_{v(t)}(r)$ and $y_{v(t)}(r)^{1/\sigma} \geq c_2(r)/n(r)$ by quasiconvexity, where $c_1, c_2$ are both slow functions. By Lemma 2.2(iv), if $\lim_{r \to 0} r^\eta \bar{n}(r)^{-1} = 0$, then $\lim_{r \to 0} r^{\eta_1} n(r)^{-1} = 0$ for all $\eta_1 > \eta$. Thus, $\beta \leq \bar{\beta}$. Let $n'(r) = \inf\{t > 0 : g_{v(t)}(r^{-1}) > m\}$. Since $X_t$ does not have the drift initially, that is, $\sum_{i=1}^d \gamma_0^{(i)}(\varepsilon) = 0$ for some $\varepsilon > 0$, $y_t(r) \approx G_t(r) + M_t(r) \approx g_t(r^{-1})$. Write $c_1 y_t(r) \leq g_t(r^{-1}) \leq c_2 y_t(r)$ where $c_1 < 1, c_2 > 1$. Fix $r$ and let $y^{-1}$ be the inverse of $y_{v(t)}(r)$. Then $n(r) = y^{-1}(m)$ and $y^{-1}(c_2^{-1}m) \leq n'(r) \leq y^{-1}(c_1^{-1}m)$. Replacing $t_i$ by $y^{-1}(t_i)$ in (3.1) yields $y^{-1}(t_2)/y^{-1}(t_1) \leq c(r)^{-1/\sigma}(t_2/t_1)^{1/\sigma}$. Applying this inequality, we find that $(c_2^{-1}c(r))^{1/\sigma} n(r) \leq n'(r) \leq (c_1 c(r))^{-1/\sigma} n(r)$. It follows from $n'(r) = \hat{n}(r^{-1})$ that $\delta = \sup\{\eta \geq 0 : \lim_{r \to \infty} r^\eta \hat{n}(r) = 0\}$, $\beta = \inf\{\eta \geq 0 : \lim_{r \to \infty} r^\eta \hat{n}(r) = \infty\}$. $\square$

Equation (3.1) is equivalent to $y_{v(t_2)}(r)/t_2 \geq c(r) y_{v(t_1)}(r)/t_1$ with $v(t)$ replaced by $v(t^{1/\sigma})$. [If $\phi(t) = t$, $y_{v(t)}(b)/t = 1$.] $y_{v(t)}(r)$ is nearly convex in $t$ since $y_{v(t_2)}(r)/t_2 \geq \rho_r y_{v(t_1)}(r)/t_1$. Chances are $\rho_r$ will drop too fast as $r$ approaches 0. ($\rho_r$ depends on $v$.) $X_t$ is said to be of class $\mathcal{I}$:

if for some $v, y_{v(t)}(r)$ is convex in $t$ for all $r$ small,

that is, $c(r) = \sigma = 1$ in (3.1).

(So, $\beta \leq 2$.) Clearly, $X_t$ is of class $\mathcal{I}$ if and only if there exist functions $h_s(r)$ nondecreasing in $s$ and $u(s)$ nondecreasing continuous with $u(0) = 0$ such that $y_t(r) = \int_0^t h_s(r) u(ds)$. In that case $v = u^{-1}$, $y_{v(t)}(r) = \int_0^t h_{v(s)}(r) \, ds$ and $u(t) = \int_0^t h_s(b)^{-1} y(b)(ds)$ for all $b$ small. Here $\int_0^t y(b)(ds) = y_t(b)$.

Class $\mathcal{I}$ is very large. Let $y_t'(r) = \frac{d}{dt} y_t(r)$. One of the conditions that $y_t(r)$ is differentiable in $t$ a.e. is that $B_t, Q_t$ (or $C_t$), $\nu_t$ each are absolutely continuous. $X_t$ is of class $\mathcal{I}$ if and only if there exists a function $g$ such that $g(s) y_s'(r)$ is nondecreasing in $s$ a.e. for all $r$, in which case $u(t) = \int_0^t g(s)^{-1} \, ds$ and $y_t(r) = \int_0^t g(s) y_s'(r) u(ds)$. Let $X_t$ be a continuous process with additive components. Then (1.3) holds. The function $u$ in the general case will be given in Section 5. For $X_t$, $y_t(r) = B_t^* r^{-1} + C_t r^{-2}$ where $B_t^* = \sum_{i=1}^d (\max_{0 \leq s \leq t} |B_s^{(i)}|)$ and $C_t = \sum_{i=1}^d C_t^{(i)}$. Assume that $\frac{d}{dt} B_t^* > 0, \frac{d}{dt} C_t > 0$ exist. $X_t$ is of class $\mathcal{I}$ if $\frac{d}{dt} C_t / \frac{d}{dt} B_t^*$ or $\frac{d}{dt} B_t^* / \frac{d}{dt} C_t$ is nondecreasing. $u(t) = B_t^*$ in the first case while $u(t) = C_t$ in the second. If $B_t^* \approx C_t$ or if $B_t^* \equiv 0$, $y_t(r) \approx C_t r^{-2}$ for $r \in (0, 1)$, in which case (3.1) holds vacuously. Hence, $\lim_{t\to 0} C_t^{-\alpha} X_t^* = 0$ or $\infty$ a.s. according as $\alpha < 1/2$ or $\alpha > 1/2$.

Let $X_t^3 = X_{f_1(t)}^1 + X_{f_2(t)}^2$ where $X_t^1, X_t^2$ are independent Lévy processes in $\mathbb{R}^d$ and $f_1, f_2$ are nondecreasing continuous functions with $f_1(0) = f_2(0) = 0$. (Locally, every additive process can be characterized as $X_t^3$.) (a) Suppose for each $r$ small, in vector terms either $m_1(r) \geq 0, m_2(r) \geq 0$ or $m_1(r) \leq 0, m_2(r) \leq 0$, where $m_i(r) = r^{-1}(B_i - \int_{r<|x|\leq 1} x\nu_i(dx)), r \in (0,1), i = 1, 2$. Then $y_t(r) = f_1(t)h_1(r) + f_2(t)h_2(r)$ for $X_t^3$ which is of class $\mathcal{I}$ if $f_1'(t) > 0, f_2'(t) > 0$ exist and one of the quotients $f_1'(t)/f_2'(t), f_2'(t)/f_1'(t)$ is nondecreasing. (b) If $f_1 \approx f_2$, then $y_t(r) \approx f_1(t)z(r)$ for $X_t^3$ and (3.1) holds.

An additive process $(X_t; B_t, Q_t, \nu_t)$ in $\mathbb{R}^d$ is a semimartingale if and only if $B_t$ is of bounded variation. Fix a nondecreasing continuous function $u$ with $u(0) = 0$, a Lévy kernel $\kappa_s(dx)$, and two $\mathbb{R}^d$-valued functions $b_s, \sigma(s)$, where $\kappa_s(dx), b_s, \sigma(s)$ are locally bounded left-continuous each. Define

$$\nu(ds, dx) = \kappa_s(dx)u(ds),$$

(3.2)

$$B_t = \int_0^t b_s u(ds), \qquad C_t = \int_0^t \sigma(s)^2 u(ds),$$

where $\sigma(s)^2 = (\sigma_1(s)^2, \sigma_2(s)^2, \ldots, \sigma_d(s)^2)$. Choose any quadratic covariation matrix $Q_t$ with $C_t$ as diagonal for a $d$-dimensional continuous Gaussian martingale. One can verify that $\nu_t(A) = \nu([0,t] \times A) = \int_0^t \kappa_s(A)u(ds)$ is a nondecreasing continuous Lévy kernel and that $B_t, C_t$ are continuous of bounded variation. Thus, a semimartingale additive process $(X_t; B_t, Q_t, \nu_t)$ is defined. Here $M_t^{(i)}(r) = |\int_0^t \tilde{m}_s^{(i)}(r)u(ds)|$ where

$$\tilde{m}_s(r) = r^{-1}\left(b_s - \int_{r<|x|\leq 1} x\kappa_s(dx)\right), \qquad s \geq 0, r < 1,$$

with components $\tilde{m}_s^{(i)}(r), 1 \leq i \leq d$ and $G_t^{(i)}(r) + K_t^{(i)}(r) = \int_0^t (r^{-2}\sigma_i(s)^2 + \int (x_i/r)^2 \wedge 1 \kappa_s^{(i)}(dx_i))u(ds)$, where $\kappa_s^{(i)}(B) = \kappa_s(\{x \in \mathbb{R}^d : x_i \in B\}), B \in \mathcal{B}(\mathbb{R})$. Let

$$h_s^{(i)}(r) = r^{-2}\sigma_i(s)^2 + \int (x_i/r)^2 \wedge 1 \kappa_s^{(i)}(dx_i) + |\tilde{m}_s^{(i)}(r)|$$

and $\mu_u$ the Lebesgue–Stieltjes measure induced by $u$. Suppose that for every $r \in (0, r_0], s \in [0, t_0], \mu_u$-a.e. with $t_0, r_0$ small, the two conditions hold:

(i) each $\tilde{m}_s^{(i)}(r)$ has no sign change in $s$;
(ii) $\sum_{i=1}^d h_s^{(i)}(r)$ is nondecreasing in $s$.

Condition (i) implies that $|\int_0^t \tilde{m}_s^{(i)}(r)u(ds)| = \int_0^t |\tilde{m}_s^{(i)}(r)|u(ds)$. Consequently, $M_t^{(i)}(r) = M_t^{(i)*}(r)$ [i.e., $M_t^{(i)}(r)$ is nondecreasing in $t$] and $y_t^{(i)}(r) =$

$\int_0^t h_s^{(i)}(r)u(ds)$. $X_t$ is of class $\mathcal{I}$ since $y_t(r) = \int_0^t (\sum_{i=1}^d h_s^{(i)}(r))u(ds)$ with $\sum_{i=1}^d h_s^{(i)}(r)$ nondecreasing in $s$ thanks to condition (ii).

Semimartingale additive processes can only be defined in that way. Given any semimartingale additive process $(X_t; B_t, Q_t, \nu_t)$ in $\mathbb{R}^d$, there exist a nondecreasing continuous function $u$ with $u(0) = 0$, a Lévy kernel $\kappa_s(dx)$, a vector $b_s$ and a nonnegative definite symmetric $d \times d$ matrix $(c_{ij}(s))$, all of which are locally bounded left-continuous, such that (3.2) holds with $q_{ij}(t) = \int_0^t c_{ij}(s)u(ds)$, that is, $c_{ii} = \sigma_i^2$. This property is better known as disintegration. By the Radon–Nikodym theorem, one can take $u(t) = \sum_{i \leq d} V_0^t B_i + \sum_{i,j \leq d} V_0^t q_{ij} + \int |x|^2 \wedge 1 \nu_t(dx), V_0^t f$ denoting the total variation of $f$ over $[0, t]$. $\kappa_s, b_s, c_{ij}(s), u(s)$ are not unique. For example, if $u$ is absolutely continuous [which implies that $y_t'(r)$ exists], any absolutely continuous nondecreasing function $u_1$ with $u_1(0) = 0$ can replace $u$ since $u(ds) = (du/ds)(du_1/ds)^{-1}u_1(ds)$. For an extensive account on the general semimartingale case, see [3], Proposition 2.9, Chapter II, page 77. The same holds true for processes with semimartingale additive components for which disintegration holds as well, that is, there exists a common function $u$ for all components. The components form an additive process in $\mathbb{R}^d$ with components independent of one another. Disintegration gives the representation of the characteristics of each component with a common function $u$. In special cases when $u_1, u_2, \ldots, u_d$ are all absolutely continuous, the $u_i$'s can be replaced by a single absolutely continuous nondecreasing function $u$ with $u(0) = 0$.

Given a semimartingale additive process $X_t$ with $M_t^{(i)}(r)$ nondecreasing in $t$, one way to argue that $X_t$ is of class $\mathcal{I}$ is to look for a combination of $\kappa_s, b_s, c_{ij}(s)$ satisfying condition (ii). [Condition (i) is equivalent to $M_t^{(i)}(r) = M_t^{(i)*}(r)$.] Note that the proof of disintegration does not include the techniques to do that. Let $X_t$ be a Lévy process with characteristics $(B, Q, \nu)$. There are many ways to represent $(B, Q, \nu)$ as (3.2); here, for instance, $\kappa_s(dx) = \nu(dx), b_s = B, \sigma_i(s)^2 = q_{ii}, u(t) = t$, or $\kappa_s(dx) = 2\sqrt{s}\nu(dx), b_s = 2\sqrt{s}B, \sigma_i(s)^2 = 2\sqrt{s}q_{ii}, u(t) = \sqrt{t}$. Clearly, conditions (i), (ii) hold.

We provide examples of $\kappa_s, b_s, c_{ij}(s)$ typically satisfying conditions (i), (ii). Let $Y_t$ be a rcll process with independent increments in $\mathbb{R}^d$ not necessarily continuous in probability. Then the Lévy kernel $\kappa_s(dx)$ induced by $Y_t$ is nondecreasing left-continuous. Take a vector $\sigma(s) = (\sigma_1(s), \sigma_2(s), \ldots, \sigma_d(s))$ with $|\sigma_i(s)|$ each nondecreasing left-continuous and let $c_{ii} = \sigma_i^2$. Conditions (i), (ii) follow if for each $r$, $\tilde{m}_s^{(i)}(r)$ is a nonnegative nondecreasing or nonpositive nonincreasing function of $s \in (0, t_0]$. To make $b_s$ match up $\int_{r < |x| \leq 1} x \kappa_s(dx)$ for all $r$, we need $\kappa_s(dx)$ to be more specific. There are a great many examples where $b_s$ can be determined but a big majority of them are rather complicated and highly irregular except for the three as follows. (Their sums are also tractable. We will see that in a moment.) (a) $\kappa_s^{(i)}$ is symmetric for every $s \in (0, t_0]$. (b) $\kappa_s^{(i)}$ is concentrated on $(0, \infty)$ or $(-\infty, 0)$ for every $s \in (0, t_0]$. (Since $\kappa_s^{(i)}$ is nondecreasing, either $(0, \infty)$ or $(-\infty, 0)$ must be held fixed for all $s \in (0, t_0]$.) (c) $\kappa_s^{(i)} = f^{(i)}(s)\nu^{(i)}, s \in (0, t_0]$ where $f^{(i)}$ is a nondecreasing function and $\nu^{(i)}$ is a Lévy measure. In case (a), $\tilde{m}_s^{(i)}(r) = r^{-1}b_s^{(i)}$. Thus, $b_s^{(i)}$ can be any nonnegative nondecreasing or nonpositive nonincreasing function. In case (b), $\tilde{m}_s^{(i)}(r) = r^{-1}(b_s^{(i)} - \int_r^1 x\kappa_s^{(i)}(dx))$, where we assume that $\kappa_s^{(i)}$ is concentrated on $(0, \infty)$. Thus, $b_s^{(i)}$ can be any nonpositive nonincreasing function. We can also take $b_s^{(i)} = \int_0^1 x\kappa_s^{(i)}(dx) + a_s^{(i)}$ where $a_s^{(i)}$ is a positive nondecreasing function if $\int_0^1 x\kappa_s^{(i)}(dx) < \infty$. In case (c), $\tilde{m}_s^{(i)}(r) = r^{-1}(b_s^{(i)} - f^{(i)}(s)\int_{r<|x|\leq 1} x\nu^{(i)}(dx))$ and we let $b_s^{(i)} = b_i f^{(i)}(s)$ where $b_i \in \mathbb{R}$.

Alternatively, we can take $\kappa_s = \nu_s, b_s = B_s, (c_{ij}(s)) = Q_s$ where $X_t$ is any additive process in $\mathbb{R}^d$ with characteristics $(B_t, Q_t, \nu_t)$ for which each $M_t^{(i)}(r)$ is nondecreasing in $t \in [0, t_0]$ for every $r \in (0, r_0]$ $[r_0 \in (0, 1)]$. [Consider $M_s^{(i)}(r)$ as $|\tilde{m}_s^{(i)}(r)|$.] $\hat{X}_t$ with $\hat{\nu}(ds, dx) = \nu_s(dx)u(ds), \hat{B}_t = \int_0^t B_s u(ds), \hat{Q}(t) = \int_0^t Q_s u(ds)$ as in (3.2) is of class $\mathcal{I}$. Since for $\hat{X}_t$, $\hat{M}_t^{(i)}(r)$ remains nondecreasing in $t$, we can take $\kappa_s = \hat{\nu}_s, b_s = \hat{B}_s, (c_{ij}(s)) = \hat{Q}_s$ and obtain another process of class $\mathcal{I}$ with the same or new function $u$.

$M_t^{(i)}(r)$ is nondecreasing in $t$ in the three cases: (a) $X_t^{(i)}$ is symmetric on $[0, t_0]$. $(M_t^{(i)}(r) \equiv 0, t \in [0, t_0], r > 0.)$ (b) $X_t^{(i)}$ is monotone on $[0, t_0]$. $[M_t^{(i)}(r) = r^{-1}(B_t^{(i)} - \int_r^1 x\nu_t^{(i)}(dx)), r \in (0, 1)$, is nondecreasing in $t$ when $X_t^{(i)}$ is increasing.] (c) $X_t^{(i)} = \bar{X}_{f(t)}^{(i)}$ where $(\bar{X}_t^{(i)}, B^{(i)}, C^{(i)}, \nu^{(i)})$ is a Lévy process and $f$ is a nondecreasing continuous function with $f(0) = 0$. $[M_t^{(i)}(r) = f(t)r^{-1}|B^{(i)} - \int_{r<|x|\leq 1} x\nu^{(i)}(dx)|.]$ Given additive processes $(X_t^i; B_t^i, C_t^i, \nu_t^i), i = 1, 2, 3$ in $\mathbb{R}$, where $X_t^1, X_t^2$ are independent and $X_t^3 = X_t^1 + X_t^2$, let $m_t^i(r) = r^{-1}(B_t^i - \int_{r<|x|\leq 1} x\nu_t^i(dx)), r \in (0, 1)$. Then $M_t^i(r) = |m_t^i(r)|$ [the function $M_t(r)$ in (2.2)], $m_t^3(r) = m_t^1(r) + m_t^2(r), B_t^3 = B_t^1 + B_t^2, C_t^3 = C_t^1 + C_t^2, \nu_t^3 = \nu_t^1 + \nu_t^2$. If $m_t^1(r), m_t^2(r)$ are both nondecreasing or both nonincreasing in $t$ for each $r$, then $M_t^3(r) = M_t^1(r) + M_t^2(r)$ and $M_t^3(r)$ is nondecreasing in $t$. For example, the following processes have gone beyond aforementioned three types: (a)+(b), (a)+(c), (a)+(b)+(c), (b)+(c), (c)+(c), and so on. If $X_t^2$ is (c), $m_t^2(r) = f(t)r^{-1}(B - \int_{r<|x|\leq 1} x\nu(dx)) = f(t)m(r)$. In (b)+(c), if $m(r) \geq 0$ for all small $r$, $X_t^1$ needs to be increasing while if $m(r) \leq 0$ for all small $r$, $X_t^1$ needs to be decreasing. In (c)+(c), $m_1(r)$ and $m_2(r)$ have to have identical signs for each $r$.

$M_t^{(i)}(r)$ is nondecreasing in $t$ in many unfamiliar cases, each of which shows its own way $\nu_t^{(i)}$ charges the area $\{x : r < |x| \leq 1\}$ in order to match up the swing of $B_t^{(i)}$. Since $M_t^{(i)}(r)$ is continuous in $t$, we can expect for a small interval $[0, t_0]$, $M_t^{(i)}(r)$ is nondecreasing in $t$ [if $\frac{d}{dt} M_t^{(i)}(r)$ is continuous, say] but the interval depends on $r$. Here is a more general case. Let $\kappa_t(dx)$ be a Lévy kernel (locally bounded left-continuous) in $\mathbb{R}$. Assume that there is $r_0 \in (0, 1)$ such that $\int_{r < |x| \leq 1} x \kappa_s(dx) \geq \int_{r_0 < |x| \leq 1} x \kappa_s(dx), r \in (0, r_0]$, for every $s \in [0, t_0]$. For a well-behaved Lévy measure $\kappa$, very often $\int_{r < |x| \leq 1} x \kappa(dx)$ shows a tendency to increase as $r \downarrow 0$. Let $b_s$ be any locally bounded left-continuous function satisfying $b_s - \int_{r_0 < |x| \leq 1} x \kappa_s(dx) \leq 0$ for every $s \in [0, t_0]$. It follows that $\tilde{m}_s(r) = r^{-1}(b_s - \int_{r < |x| \leq 1} x \kappa_s(dx)) \leq 0$ for every $r \in (0, r_0], s \in [0, t_0]$, which is condition (i), and that $M_t(r) = -\int_0^t \tilde{m}_s(r) u(ds)$ is nondecreasing in $t \in [0, t_0]$ for all $r \in (0, r_0)$.

We also offer an example where $y_t(r) \approx f(t) z(r)$. Assume $d = 1$ first. Since condition (i) implies that $y_t(r) = \int_0^t h_s(r) u(ds), r \in (0, r_0], t \in [0, t_0]$, if $h_s(r) \approx c_s z(r)$, $y_t(r) \approx f(t) z(r)$ where $f(t) = \int_0^t c_s u(ds)$. Let $(b_s, \sigma(s)^2, \kappa_s)$ be the characteristics of a centered Lévy process in $L^2$; that is, $\int_{|x| > 1} x^2 \kappa_s(dx) < \infty$ and $b_s = -\int_{|x| > 1} x \kappa_s(dx)$. Then $\tilde{m}_s(r) = r^{-1}(b_s - \int_{r < |x| \leq 1} x \kappa_s(dx)) = -r^{-1} \int_{|x| > r} x \kappa_s(dx)$ and $h_s(r) < 2(\sigma(s)^2 + \int x^2 \kappa_s(dx)) r^{-2}$. Thus, $h_s(r) \approx \sigma(s)^2 r^{-2}$ if $\int x^2 \kappa_s(dx) \leq \hat{c} \sigma(s)^2$. Define $\kappa_s(dx) = c_s'(2 - \alpha_s) x^{-(1+\alpha_s)} dx, x \in (0, 1], \kappa_s(dx) = c_s''(2 - \alpha_s) |x|^{-(1+\alpha_s)} dx, x \in [-1, 0)$, and $\kappa_s = 0$ on $\{x : |x| > 1\}$ where $\alpha_s, c_s', c_s''$ are continuous functions on $[0, t_0]$ satisfying $\alpha_s \in [0, 2], c_s' > c_s'' > 0$. Then $\int_{|x| > r} x \kappa_s(dx) = \int_{r < |x| \leq 1} x \kappa_s(dx) = (c_s' - c_s'') B(r, \alpha_s) \geq 0$ (nonincreasing in $r$) where $B(r, \alpha_s) = (2 - \alpha_s) |\alpha_s - 1|^{-1} |r^{1-\alpha_s} - 1|$ or $B(r, \alpha_s) = \log(1/r)$ according as $\alpha_s \neq 1$ or $\alpha_s = 1$. Thus, for each $r$, $\tilde{m}_s(r) \leq 0$ as $s$ varies. Let $b_s = 0$ and $\sigma(s)$ a continuous function on $[0, t_0]$ satisfying $\sigma(s)^2 \geq \hat{c}(c_s' - c_s'') = \hat{c} \int x^2 \kappa_s(dx)$ for some fixed constant $\hat{c} \in (0, \infty)$. The reader can verify that $\nu_t([r, 1]) = \nu([0, t] \times [r, 1]) = \int_0^t \kappa_s([r, 1]) u(ds)$ is continuous in $t$ for every fixed $r \in (0, 1)$. Thus, $\nu(ds, dx) = \kappa_s(dx) u(ds), B_t = 0, C_t = \int_0^t \sigma(s)^2 u(ds)$ determine an additive process $X_t$ (a martingale with the jump size bounded by 1) for which $y_t(r) \approx C_t r^{-2}, r \in (0, 1), t \in [0, t_0]$. Clearly $\delta = \beta = 2$ in this example. One can easily make a similar example in the case $d > 1$ and even an example where $X_t$ has large jumps. The interested reader would probably demand a more interesting example where $\int_{|x| > r} x \kappa_s(dx)$ takes both positive and negative values as $r$ varies, Gaussian part is not so prominent, $z(r)$ is sophisticated enough to force $\delta < \beta$, and so on, or even $M_t(r)$ fails to be nondecreasing in $t$. Unfortunately, the new technique for that is at present unavailable.

Let $e^{-1}$ be the inverse of a continuous moderate function $e$. If $y_{v(t)}(r)$ is quasiconvex, so too is $y_{v \circ e^{-1}(t)}(r)$. Thus, Theorem 3.2 holds for all functions in the form $e \circ u$. Same goes here: If $y_{v(t)}(r)$ is not quasiconvex for $v$ satisfying $y_{v(t)}(b) = t$, that is, $v$ is the inverse of $y_t(b)$, then $y_{v(t)}(r)$ cannot be quasiconvex for any $v$ satisfying $y_{v(t)}(b) = \phi(t)$ where $\phi$ is moderate. $[t^p(\log(1/t))^\kappa, t^p(\log\log(1/t))^\kappa, p > 0, \kappa \in \mathbb{R}$, etc. are both quasiconvex and moderate but it is not so easy to give an example of nonmoderate nonexponential-type quasiconvex functions. If $\phi$ is exponential, $\delta_1 = \delta_2 = \beta_1 = \beta_2 = 0$ since $\sigma$ can be arbitrarily large.] In the case that $y_t(r) \approx f(t)z(r)$ where $f$ is some nondecreasing continuous function with $f(0) = 0$ and $z$ is a positive function, if we let $v = f^{-1}$ the inverse of $f$ $(u = f)$, then $y_{v(t_2)}(r)/y_{v(t_1)}(r) \approx t_2/t_1$. Equation (3.1) holds tautologically. Thus, in Theorem 3.2, $u = e \circ f$ with $\delta = \inf\{\eta \geq 0 : \liminf_{r \to 0} r^\eta/e(z(r)^{-1}) = 0\}, \beta = \inf\{\eta \geq 0 : \lim_{r \to 0} r^\eta/e(z(r)^{-1}) = 0\} \leq 2\sigma$ where $\sigma$ is an exponent for $e$. $(\delta, \beta$ do not depend on $f$.) Particularly, for a Lévy process $X_t$, Pruitt's result is extended from the case $t^{-1/\eta}X_t^*$ to the general case $e(t)^{-1/\eta}X_t^*$ with $\delta = \inf\{\eta \geq 0 : \liminf_{r \to 0} r^\eta/e(h(r)^{-1}) = 0\}, \beta = \inf\{\eta \geq 0 : \lim_{r \to 0} r^\eta/e(h(r)^{-1}) = 0\}$. A lower function for a Lévy process $X_t$ is a moderate function $e$ satisfying $\liminf_{t \to 0} e(t)^{-1}X_t^* = c \in (0, \infty)$ a.s. Assume that $X_t$ is not a compound Poisson process. Then $h(r)^{-1} \approx \mathsf{k}(r) = r^{-1}\int_0^r \sup_{0 < s \leq t} h(s)^{-1}\,dt$ which is a strictly increasing absolutely continuous moderate function with $\mathsf{k}(\infty) = \infty$. If the inverse $\mathsf{k}^{-1}$ is moderate, $\delta = \beta = 1$ with $e = \mathsf{k}^{-1}$, which implies that $\lim_{t \to 0} \mathsf{k}^{-1}(t)^{-\alpha}X_t^* = 0$ or $\infty$ a.s. according as $\alpha < 1$ or $\alpha > 1$. We suspect the lower function exists only when $\mathsf{k}^{-1}$ is moderate, in which case $e(t) = \mathsf{k}^{-1}(t)g(t)$ is a lower function where $g$ is moderate satisfying $\lim_{t \to 0} \mathsf{k}^{-1}(t)^p g(t)^q = 0$ for all $p > 0, q \in \mathbb{R}$.

There are also results for $t \to \infty$ analogous to Theorems 3.1, 3.2 as long as $\lim_{t \to \infty} X_t^* = \infty$ a.s., for if $T_r = \inf\{t > 0 : |X_t| > r\}, r > 0$, is infinite for large $r$, the probability that $\limsup_{t \to \infty} t^{-1/\eta}X_{v(t)}^* = \infty$ will be less than 1 for any function $v$ and power $\eta$. For additive processes, that can happen. For example, $P(T_r = \infty) > 0$, $r$ being large for the process $X_t^f = X_{f(t)}$ with $f$ bounded and in the case $y_t(r) \approx f(t)z(r)$, if $f$ is bounded, $P(T_r < \infty) < 1$ for large $r$. Technically, we need to reverse the symbols used for $t \to 0$ to get the results for $t \to \infty$, including such changes as $t \to 0$ to $t \to \infty$, $r \to 0$ to $r \to \infty$, ">" to "<" and vice versa, and "inf" to "sup" and vice versa. Accordingly, a sequence $\sigma_n \uparrow \infty$ is called the $\Sigma$-sequence if $\frac{\sigma_{n-1}}{\sigma_n} \cdot \sigma_n^\eta \to \infty$ as $n \to \infty$ for all $\eta > 0$ and $y_{v(t)}(r)$ is called quasiconvex if (3.1) holds for $t_0 \leq t_1 < t_2, r_0 \leq r$ with $c(r)$ replaced by $c(r^{-1})$. For Theorem 3.1, we assume that there is a sequence $v_n \uparrow \infty$ and some $\kappa \in (0, \infty)$ such that $\sum y_{v_n}(\bar{\sigma}_n^{1/\kappa})^{-1} < \infty$ (which guarantees that $T_r < \infty$ for all $r > 0$). That holds if and only if $y_\infty(r_1) = \infty$ for some $r_1$, which is equivalent to $y_\infty(r) = \infty$ for all $r$ by (2.5). [That also implies that $v$ satisfying $y_{v(t)}(b) = \phi(t)$ is finitely determined for all quasiconvex functions $\phi$ with $\phi(\infty) = \infty$.] The order of the indices is reversed as $\delta_1 \leq \delta_2 \leq 2/\sigma, \beta_2 \leq \beta_1, \beta_2 \leq \delta_1, \beta_1 \leq \delta_2$. As far as Theorem 3.2 goes, we assume that $y_{v(t)}(r)$ is quasiconvex in the $t \to \infty$ sense with $v(\infty) = \infty$. Analogously, $\delta = \delta_P = \delta_E$, $\beta = \beta_P$ with $\int_0^l$ replaced by $\int_l^\infty$ in addition.

**4. The semicontinuous function method.** For $b \in (0, \infty)$ fixed, define

$$\delta = \inf\left\{\eta \geq 0 : \liminf_{r \to 0} r^\eta y_b(r) = 0\right\}, \qquad \beta = \inf\left\{\eta \geq 0 : \lim_{r \to 0} r^\eta y_b(r) = 0\right\}.$$

By (2.5), $y_b(r) \leq 3r^{-2} y_b(1)$ for $r \in (0, 1)$. Thus, $\delta \leq \beta \leq 2$. Define for $t$ small,

$$\overline{v}(t) = \inf\{s > 0 : y_s(t^{1/\eta}) \geq c(t) t y_b(t^{1/\eta}) \text{ for all } \eta \in [\delta - \varepsilon, \delta)\},$$
$$\underline{v}(t) = \sup\{s > 0 : y_s(t^{1/\eta}) \leq c(t)^{-1} t y_b(t^{1/\eta}) \text{ for all } \eta \in (\delta, \delta + \varepsilon]\},$$
$$\overline{u}(t) = \inf\{s > 0 : y_s(t^{1/\eta}) \geq c(t) t y_b(t^{1/\eta}) \text{ for all } \eta \in [\beta - \varepsilon, \beta)\},$$
$$\underline{u}(t) = \sup\{s > 0 : y_s(t^{1/\eta}) \leq c(t)^{-1} t y_b(t^{1/\eta}) \text{ for all } \eta \in (\beta, \beta + \varepsilon]\},$$

where $\varepsilon$ is a small positive constant and $c(t)$ is a continuous slow function. $\overline{v}, \underline{v}, \overline{u}, \underline{u}$ are finitely determined positive (but not necessarily monotone) functions.

THEOREM 4.1. *Let $X_t$ be any process with additive components and $\delta$, $\beta$, $\overline{v}$, $\underline{v}$, $\overline{u}$, $\underline{u}$ as given above. Then with probability* 1 (1.2a), (1.2b) *hold.*

PROOF. Let $\hat{v}(t) = \inf\{s > 0 : \dot{y}_s(t^{1/\eta}) \geq c(t) t \dot{y}_b(t^{1/\eta}) \text{ for all } \eta \in [\delta - \varepsilon, \delta)\}$. Suppose that $\lim_{n \to \infty} \hat{v}(s_n) = x$ for a sequence $s_n \to t$. Since $\dot{y}_{\hat{v}(s_n)}(s_n^{1/\eta}) \geq c(s_n) s_n \dot{y}_b(s_n^{1/\eta})$ and since $\dot{y}_t(r)$ is jointly continuous by Lemma 2.3 and $c(t)$ is continuous, $\dot{y}_x(t^{1/\eta}) \geq c(t) t \dot{y}_b(t^{1/\eta})$. Thus, $x \geq \hat{v}(t)$, which shows that $\hat{v}$ is lower semicontinuous. Since $y_s(t^{1/\eta}) \geq c(t) t y_b(t^{1/\eta})$ implies $\dot{y}_s(t^{1/\eta}) \geq kc(t) t \dot{y}_b(t^{1/\eta})$ where $k \in (0, \infty)$ is a constant, $\overline{v} \geq \hat{v}$. [The slow function for $\hat{v}$ is $kc(t)$ not $c(t)$ now.] Since $\hat{v}$ is lower semicontinuous, there is $\sigma_n \in [2^{-(n+1)}, 2^{-n}]$ such that $\hat{v}(\sigma_n) \leq \hat{v}(t)$ for all $t \in [2^{-(n+1)}, 2^{-n}]$. Suppose that $\eta < \eta_1 < \eta_2 < \delta$ with $\eta_1 \geq \delta - \varepsilon$. Since $\eta_2 < \delta$, $y_b(r^{1/\eta_1})^{-1} \leq c_1 r^{\eta_2/\eta_1}$ for all small $r$. Since $\dot{y}_{\hat{v}(t)}(t^{1/\eta}) \geq c(t) t \dot{y}_b(t^{1/\eta})$ for all $\eta \in [\delta - \varepsilon, \delta)$ and $y_t(r) \approx \dot{y}_t(r)$, by the second bound in (2.6) with $k = 2$, for small $t$, $P(X^*_{\hat{v}(t)} \leq t^{1/\eta_1}) \leq c_2(y_{\hat{v}(t)}(t^{1/\eta_1}))^{-1} \leq c_3 c(t)^{-1} (t y_b(t^{1/\eta_1}))^{-1} \leq c_4 c(t)^{-1} t^{\eta_2/\eta_1 - 1} \leq c_4 t^{\eta_2/\eta_1 - 1 - \theta}$ with $\eta_2/\eta_1 - 1 - \theta > 0$. Hence, $\sum P(X^*_{\hat{v}(\sigma_n)} \leq \sigma_n^{1/\eta_1}) < \infty$ and $X^*_{\hat{v}(\sigma_n)} > \sigma_n^{1/\eta_1}$ for large $n$ a.s. by the Borel–Cantelli lemma. Thus, for small $t$ with $t \in [2^{-(n+1)}, 2^{-n}]$, $X^*_{\hat{v}(t)} \geq X^*_{\hat{v}(\sigma_n)} > \sigma_n^{1/\eta_1} \geq 2^{-1/\eta_1} t^{1/\eta_1}$, which yields $\lim_{t \to 0} t^{-1/\eta} X^*_{\hat{v}(t)} = \infty$ a.s. We obtain the first half of (1.2a) since $\overline{v} \geq \hat{v}$. For the latter half of (1.2a), let $\check{v}(t) = \sup\{s > 0 : \dot{y}_s(t^{1/\eta}) \leq c(t)^{-1} t \dot{y}_b(t^{1/\eta}) \text{ for all } \eta \in (\delta, \delta + \varepsilon]\}$. Note that $\check{v}$ is upper semicontinuous and $\underline{v} \leq \check{v}$. The rest follows the first lead as always. □

Very often $\overline{v}/\underline{v} \leq 1, \overline{u}/\underline{u} \leq 1$ over a small interval $(0, t_0)$, which will be justified below. If that is the case, (1.3) holds for $v_i, v_s$ satisfying $\overline{v} \leq v_i \leq$ $\underline{v}, \overline{u} \le v_s \le \underline{u}$. $\overline{v}(t)/\underline{v}(t), \overline{u}(t)/\underline{u}(t)$ are bounded by $t^{-\alpha}$ where $\alpha \downarrow 0$ as $\varepsilon \downarrow 0$. Let $\overline{v}_\theta(t) = \inf\{s > 0 : y_s(t^{1/\delta}) \ge c(t)t^{\theta/\delta}y_b(t^{1/\delta})\}(\theta < \delta)$, $\underline{v}_\theta(t) = \sup\{s > 0 : y_s(t^{1/\delta}) \le c(t)^{-1}t^{\theta/\delta}y_b(t^{1/\delta})\}(\theta > \delta)$, $\overline{u}_\theta(t) = \inf\{s > 0 : y_s(t^{1/\beta}) \ge c(t)t^{\theta/\beta}y_b(t^{1/\beta})\}(\theta < \beta)$ and $\underline{u}_\theta(t) = \sup\{s > 0 : y_s(t^{1/\beta}) \le c(t)^{-1}t^{\theta/\beta}y_b(t^{1/\beta})\}(\theta > \beta)$. If $s$ satisfies $y_s(t^{1/\delta}) \ge c(t)t^{\theta/\delta}y_b(t^{1/\delta})$ for $\theta < \delta, t < 1$, then by (2.5) for all $\eta \in [\delta - \varepsilon, \delta)$ and $\varepsilon$ small enough, $y_s(t^{1/\eta}) \ge kc(t)t^{\theta/\delta + 2(1/\eta - 1/\delta)}y_b(t^{1/\eta}) \ge kc(t)ty_b(t^{1/\eta})$ where $k \in (0,1)$ is some constant. Thus, $b > \overline{v}_\theta \ge \overline{v} > 0$. Similarly, $b > \underline{v} \ge \underline{v}_\theta > 0$, $b > \overline{u}_\theta \ge \overline{u} > 0$, $b > \underline{u} \ge \underline{u}_\theta > 0$. Equation (1.2c) follows with $v(\eta, t) = \inf\{s > 0 : y_s(t^{1/\delta}) > t^{\eta/\delta}y_b(t^{1/\delta})\}$, $u(\eta, t) = \inf\{s > 0 : y_s(t^{1/\beta}) > t^{\eta/\beta}y_b(t^{1/\beta})\}$, $t \in (0,1)$, $\eta > 0$.

PROPOSITION 4.2. *Any of the three conditions below implies that $\overline{v} \le \underline{v}$, $\overline{u} \le \underline{u}$ in Theorem 4.1.*

(i) $y_s(r_2)/y_b(r_2) \le c_1(r_1)^{-1}c_2(r_2)^{-1}y_s(r_1)/y_b(r_1)$ for all $r_1 < r_2$ small and $s \in (0,b)$.

(ii) $y_s(r)/\dot{y}_s(r) - y_b(r)/\dot{y}_b(r) \le \log c(r)/\log r$ for all small $r$ and $s \in (0,b)$.

(iii) *For some continuous function $v > 0$ with values $v(s_n) = \inf\{s > 0 : \dot{y}_s(s_n^{1/\delta}) \ge c(s_n)s_n\dot{y}_b(s_n^{1/\delta})\}$ at a sequence $s_n \downarrow 0$, $\overline{f}(t, \eta) = [\dot{y}_{v(t)}(t^{1/\eta})/c(t)t \times \dot{y}_b(t^{1/\eta})] \wedge 1, t > 0, \eta > 0$ ($\underline{f}(t, \eta) = [c(t)t\dot{y}_b(t^{1/\eta})/\dot{y}_{v(t)}(t^{1/\eta})] \wedge 1$) is uniformly continuous on $(0, t^*] \times [\delta - \varepsilon, \delta]$ ($(0, t^*] \times [\delta, \delta + \varepsilon]$). (There is a similar statement for $\overline{u}, \underline{u}$ as well.)*

[Here $c_1(t), c_2(t), c(t)$ are continuous slow functions.]

PROOF. $y_s(r)/y_b(r)$ ($s < b$) and $y_t(r)/\dot{y}_t(r)$ are two interesting functions. Neither is monotone in $r$. The former takes values in $(0,1]$ while the latter vacillates between $48^{-1}$ and 2. Since $\dot{y}_t(r)$ is jointly continuous, $\overline{f}$ is continuous (and bounded by 1) but not necessarily uniformly continuous in the region $(0, t^*] \times [\delta - \varepsilon, \delta]$.

(i) For $t \in (0,1)$, let $s < b$ be such that $y_s(t^{1/\delta})/y_b(t^{1/\delta}) = t$. With the condition in (i) we have both $y_s(t^{1/\eta})/y_b(t^{1/\eta}) \ge c_1(t^{1/\eta_1})c_2(t^{1/\eta_1})t$ for all $\eta \in [\eta_1, \delta)(\eta_1 = \delta - \varepsilon)$ and $y_s(t^{1/\eta})/y_b(t^{1/\eta}) \le (c_1(t^{1/\eta_2})c_2(t^{1/\eta_2}))^{-1}t$ for all $\eta \in (\delta, \eta_2](\eta_2 = \delta + \varepsilon)$, where both $c_1(r^{1/\eta_1})c_2(r^{1/\eta_1})$ and $c_1(r^{1/\eta_2})c_2(r^{1/\eta_2})$ are slow functions. Thus, $\overline{v} \le \underline{v}$.

(ii) Let $g(\eta) = \log(\dot{y}_s(t^{1/\eta})/\dot{y}_b(t^{1/\eta}))$ for $\eta > 0$ with $s > 0, t > 0, b > 0$ all fixed. Thanks to $\frac{d}{dr}I_t(r) = r^{-1}(y_t(r)^{-1} - I_t(r))$ and the mean value theorem, $\dot{y}_s(t^{1/\eta})/\dot{y}_b(t^{1/\eta}) = e^{g'(\theta)(\eta - \eta')}\dot{y}_s(t^{1/\eta'})/\dot{y}_b(t^{1/\eta'})$ for $\eta < \eta'$, where $\theta \in (\eta, \eta')$, $g'(\theta) = C_\theta \theta^{-1}\log(1/t^{1/\theta})(y_s(t^{1/\theta})/\dot{y}_s(t^{1/\theta}) - y_b(t^{1/\theta})/\dot{y}_b(t^{1/\theta}))$, and $C_\theta = (\dot{y}_b(t^{1/\theta})\dot{y}_s(t^{1/\theta}))/(y_b(t^{1/\theta})y_s(t^{1/\theta})) \in [2^{-2}, 48^2]$. For $t \in (0,1)$, let $s < b$ be such that $y_s(t^{1/\delta})/y_b(t^{1/\delta}) = t$. The assumption of (ii) implies that $e^{g'(\theta)(\eta - \delta)} \ge c(t^{1/\eta_1})^{48^2\varepsilon/\eta_1}$ for all $\eta \in [\eta_1, \delta)(\eta_1 = \delta - \varepsilon)$ and $e^{g'(\theta)(\delta - \eta)} \ge c(t^{1/\delta})^{48^2\varepsilon/\delta}$ for all $\eta \in (\delta, \delta + \varepsilon]$. It follows that $y_s(t^{1/\eta})/y_b(t^{1/\eta}) \geq \rho c(t^{1/\eta_1})^{48^2\varepsilon/\eta_1} t$ for all $\eta \in [\delta - \varepsilon, \delta)$ and $t \geq \rho c(t^{1/\delta})^{48^2\varepsilon/\delta} y_s(t^{1/\eta})/y_b(t^{1/\eta})$ for all $\eta \in (\delta, \delta + \varepsilon]$, where $\rho \in (0,1)$ is a constant and $c(r^{1/\eta_1})^{48^2\varepsilon/\eta_1}, c(r^{1/\delta})^{48^2\varepsilon/\delta}$ each are slow functions. Thus, $\overline{v} \leq \underline{v}$.

(iii) The assumption here implies that $\overline{f}$ has a continuous extension to $[0, t^*] \times [\delta - \varepsilon, \delta]$. In particular, $\overline{f}$ is continuous at point $(0, \delta)$. By the definition of $v(s_n)$, $\overline{f}(s_n, \delta) = 1$. Hence, $\overline{f}(0, \delta) = 1$ and there is a neighborhood $\mathcal{O} = (0, t') \times (\delta - \varepsilon', \delta)$ of $(0, \delta)$ such that $\overline{f}(t, \eta) \geq c' > 0$, that is, $\dot{y}_{v(t)}(t^{1/\eta})/c(t)t\dot{y}_b(t^{1/\eta}) \geq 1 \wedge c'$, for all points $(t, \eta) \in \mathcal{O}$, which yields $\overline{v} \leq v$. Same goes for $\underline{f}$ and $v \leq \underline{v}$ follows. The proof that $\overline{u} \leq \underline{u}$ in each case (i), (ii), (iii) proceeds analogously. □

Here, for instance, in the case $y_t(r) \approx f(t)z(r)$ including the Lévy process case $y_t(r) = th(r)$, we can assume $y_t(r) = f(t)z(r)$. All the three conditions in Proposition 4.2 hold, where slow functions are constants in $(0,1)$. We see that $\overline{v} \leq f^{-1} \leq \underline{v}, \overline{u} \leq f^{-1} \leq \underline{u}$, where $f^{-1}$ is the inverse of $f$.

There seems no way to know the sign of $\limsup_{r \to 0} (y_s(r)/\dot{y}_s(r) - y_b(r)/\dot{y}_b(r))$. Part (ii) of Proposition 4.2 implies that $\limsup_{r \to 0}(y_s(r)/\dot{y}_s(r) - y_b(r)/\dot{y}_b(r)) \leq 0$. Replacing $y_t(r)$ by $\dot{y}_t(r)$ if necessary, we can assume $y_t(r)$ is continuous in $r$. For each $r \in (0,1)$ there is $L_r \in (0, \infty)$ such that $y_s(r)/y_s(x) - y_b(r)/y_b(x) < (\log(1/r))^{-1}$ for $x \in [r^{L_r^{-1}+1}, r)$. Since $r^p \leq -(ep)^{-1}(\log r)^{-1}$ for $p > 0, r \in (0,1)$, $y_s(r)/\dot{y}_s(r) - y_b(r)/\dot{y}_b(r) = r^{-1} \int_0^r (y_s(r)/y_s(x) - y_b(r)/y_b(x))\, dx = r^{-1} \int_0^{r^{L_r^{-1}+1}} + r^{-1} \int_{r^{L_r^{-1}+1}}^r < -(2e^{-1}L_r + 1)/\log r$. Thus, (ii) of Proposition 4.2 holds if $L_r \leq -\log c(r)$ for some slow function $c(r)$. Part (ii) of Proposition 4.2 holds on a number of occasions and yet our calculations just came up short. We are unable to pass a judgment on $L_r$.

Since $\overline{v}_\theta \geq \overline{v}, \underline{v} \geq \underline{v}_\theta$, $\overline{v}, \underline{v}$ can be replaced by $\overline{v}_\theta, \underline{v}_\theta$, respectively, in (1.2a), which remains valid as $\theta \to \delta$ while $\overline{v}_\theta \downarrow \overline{v}_\delta = \lim_{\theta \uparrow \delta} \overline{v}_\theta$, $\underline{v}_\theta \uparrow \underline{v}_\delta = \lim_{\theta \downarrow \delta} \underline{v}_\theta$, $\overline{v}_\delta \leq v \leq \underline{v}_\delta$ since $c(r) \in (0,1)$ where $v(t) = \inf\{s > 0 : y_s(t^{1/\delta}) \geq ty_b(t^{1/\delta})\}$. Thus, Theorem 4.1 is asymptotically optimal. But we are just unable to push the argument more, that is, taking limits under (1.2a) to obtain the exact result. There is another way to obtain an asymptotically optimal result. Define for $\varepsilon \in (0, 1/2)$, $\delta^\varepsilon = \inf\{\eta \geq 0 : \liminf_{r \to 0} r^{\eta(1+2\varepsilon)} y_b(r) = 0\}$, $\delta_\varepsilon = \inf\{\eta \geq 0 : \liminf_{r \to 0} r^{\eta(1-2\varepsilon)} y_b(r) = 0\}$, $v^\varepsilon(t) = \inf\{s > 0 : y_s(t^{1/\delta^\varepsilon}) > ty_b(t^{1/\delta^\varepsilon})\}$, $v_\varepsilon(t) = \inf\{s > 0 : y_s(t^{1/\delta_\varepsilon}) > ty_b(t^{1/\delta_\varepsilon})\}$. Then $\delta = \delta^\varepsilon(1 + 2\varepsilon) = \delta_\varepsilon(1 - 2\varepsilon)$. In the proof of Proposition 4.2 we derived $\dot{y}_s(t^{1/\eta})/\dot{y}_b(t^{1/\eta}) = t^{(\eta'-\eta)\theta^{-2}C_{t,\theta}}\dot{y}_s(t^{1/\eta'})/\dot{y}_b(t^{1/\eta'})$ for $\eta < \eta'$ where $\theta \in (\eta, \eta')$, $C_{t,\theta} = \dot{y}_b(t^{1/\theta})/y_b(t^{1/\theta}) - \dot{y}_s(t^{1/\theta})/y_s(t^{1/\theta})$ with $|C_{t,\theta}| < 47.5$. It follows that $\lim_{t \to 0} t^{-1/\eta} X^*_{v^\varepsilon(t)} = \infty$ a.s. if $\eta < \delta^\varepsilon$ while $\liminf_{t \to 0} t^{-1/\eta} X^*_{v_\varepsilon(t)} = 0$ a.s. if $\eta > \delta_\varepsilon$. But if $\eta < \delta$ ($\eta > \delta$), $\eta < \delta^\varepsilon < \delta$ ($\eta > \delta_\varepsilon > \delta$) for all small $\varepsilon$ while both $v^\varepsilon$ and $v_\varepsilon$ converge to $v$. The exact result follows if there is a function that can replace both $v^\varepsilon$ and $v_\varepsilon$ for a sequence $\varepsilon_n \downarrow 0$. But $v^\varepsilon, v^{\varepsilon'}, v_\varepsilon, v_{\varepsilon'}(\varepsilon \ne \varepsilon'), v, \overline{v}, \underline{v}$ just do not have a $\le$ or $\ge$ relationship between any two of them.

$\delta, \beta$ are due for simplification. By (ii), (iii) of Lemma 2.2, $\delta = \beta = 2$ if $\sum_{j=1}^d C_b^{(j)} > 0$ and $\beta \le 1$ if $\sum_{j=1}^d C_b^{(j)} = 0$ and $\sum_{j=1}^d \int_{|x| \le 1} |x| \nu_b^{(j)}(dx) < \infty$ while $\delta = \beta = 1$ if $\sum_{j=1}^d \gamma_0^{(j)*}(b) > 0$. It remains to consider initially drift-free processes $X_t$ with $\sum_{j=1}^d C_\varepsilon^{(j)} = 0$ for some $\varepsilon > 0$. By Lemma 2.2(iv), $\beta = \inf\{\eta \ge 0 : \lim_{r \to 0} r^\eta G_b(r) = 0\} = \inf\{\eta \ge 0 : \sum_{j=1}^d \int_{|x| \le 1} |x|^\eta \nu_b^{(j)}(dx) < \infty\}$. The same does not hold for $\delta$. Nonetheless under the sector condition $|\mathrm{Im}\Psi_t(\lambda)/\mathrm{Re}\Psi_t(\lambda)| \le c_0$ (which is not a sample-path type condition), $\delta$ can be defined in terms of $\mathrm{Re}\Psi_b(\lambda)$; see [7], Example 5.5(1), page 595. If $X_t$ is a process with increasing additive components and no drift on $[0, b]$, we have $\delta = \sup\{\eta \ge 0 : \lim_{r \to \infty} r^{-\eta} g_b(r) = \infty\}, \beta = \inf\{\eta \ge 0 : \lim_{r \to \infty} r^{-\eta} g_b(r) = 0\} \le 1$. To obtain the result for $t \to \infty$ analogous to Theorem 4.1, in addition to symbol reversal we need to assume that for some constant $\alpha \in (0, \infty)$ and each large $t$ there is $s \in (0, \infty)$ such that $y_s(t^{1/\delta}) \ge t^{1+\alpha} y_b(t^{1/\delta})$ and the same holds when $\delta$ is replaced by $\beta$. The assumption holds vacuously if $y_\infty(r_1) = \infty$ for some $r_1 > 0$. $\overline{v}, \underline{v}, \overline{u}, \underline{u}$ remain the same as before except for a position switch between $[\delta - \varepsilon, \delta)$ and $(\delta, \delta + \varepsilon]$ as well as between $[\beta - \varepsilon, \beta)$ and $(\beta, \beta + \varepsilon]$. $\beta, \delta$ satisfy $\beta \le \delta \le 2$ by (2.5). To simplify $\beta, \delta$ we follow the note for $r \uparrow \infty$ to Lemma 2.2. Assume that $d = 1$ for simplicity. If $X_b \in L^1$, or equivalently $\int_{|x| > 1} |x| \nu_b(dx) < \infty$, then $1 \le \beta \le \delta$ while $\beta = \delta = 1$ if $EX_s \ne 0$ for some $s \in (0, b]$ by (c). If $X_t$ is in $L^2$ and centered on $[0, b]$, or equivalently $\int_{|x| > 1} x^2 \nu_b(dx) < \infty$, and $EX_s = 0$ for all $s \in (0, b]$, then $\beta = \delta = 2$ by (b). Obviously, in any event if $C_b > 0$, $\beta = \delta = 2$. Finally, assuming $C_b = 0$, if $E|X_s| = \infty$ for all $s \in (0, b]$, or if $E|X_b| < \infty$ and $EX_s = 0$ for all $s \in (0, b]$, then $\beta = \sup\{\eta \in [0, 2] : \lim_{r \to \infty} r^\eta G_b(r) = 0\} = \sup\{\eta \in [0, 2] : \int_{|x| > 1} |x|^\eta \nu_b(dx) < \infty\}$.

**5. The moment method.** This method is not that far away from the framework of the law of the iterated logarithm for the sum of arbitrary independent r.v.'s where the growth function $u$ up to a log log term is chosen from the moments of the process. Let $X_t$ be any process in $\mathbb{R}^d$ with $X_0 = 0$ and let $e$ be a bounded nondecreasing function with $e(0) = 0$. Define

$$a_e(t) = Ee(X_t^*), \qquad H(r) = e(r)^{-1}, \qquad h(r) = (Ea_e(T_r))^{-1},$$

where $T_r = \inf\{t > 0 : |X_t| > r\}, r > 0$. By Markov's inequality, $P(X_t^* \ge r) \le a_e(t) H(r), P(X_t^* \le r) \le P(T_r \ge t) \le (a_e(t) h(r))^{-1}$. If $e$ is absolutely continuous and $X_t$ is rcll and continuous in probability, $a_e$ is a continuous function and hence $P(X_{v(t)}^* \ge r) \le tH(r), P(X_{v(t)}^* \le r) \le (th(r))^{-1}$, where $v$ is the inverse of $a_e$. If we define

$$\delta_1 = \inf\left\{\eta \ge 0 : \liminf_{r \to 0} r^\eta H(r) = 0\right\},$$

$$\delta_2 = \inf\left\{\eta \geq 0 : \liminf_{r \to 0} r^\eta h(r) = 0\right\},$$

$$\beta_1 = \inf\left\{\eta \geq 0 : \lim_{r \to 0} r^\eta h(r) = 0\right\},$$

$$\beta_2 = \inf\left\{\eta \geq 0 : \lim_{r \to 0} r^\eta H(r) = 0\right\}, \qquad u = a_e,$$

we obtain (1.1). Clearly, if $H(r) \leq c(r)^{-1} h(r)$ for some slow function $c(r)$, then $\delta_1 = \delta_2, \beta_1 = \beta_2$. In the case of additive processes, $h$ can be worked out explicitly. If $e$ is moderate, so too is $e \wedge e(1)$. The latter is bounded. On top of that, the strictly increasing absolutely continuous moderate function $\hat{e}(t) = t^{-1} \int_0^t e(x)\, dx, t > 0, \hat{e}(0) = 0$, satisfies $\hat{e} \approx e$ where the constants in $\approx$ depend only on $\rho, \sigma$. So, $e$ will be considered as bounded $[e(r) = e(1)$ for $r > 1]$ absolutely continuous for the time being. Clearly, for $e$ moderate, $\beta_2 \leq \sigma$. de la Peña and Eisenbaum [5] showed that for any rcll process $X_t$ in $\mathbb{R}^d$ with independent increments and any moderate function $e$, $Ee(X_T^*) \approx Ea_e(T)$ over *all* stopping times $T$ with the constants in $\approx$ depending on $e$ only. That being said, for stopping times $T_r$ we have $Ea_e(T_r) \approx Ee(X_{T_r}^*) = Ee(|X_{T_r}|)$. The complete result on the growth behavior of a continuous additive process $X_t$ in $\mathbb{R}^d$ is now available. Since $T_r < \infty$ a.s. for $r$ small and $|X_{T_r}| = r$, $Ee(|X_{T_r}|) = e(r)$ for any function $e$. Thus, $h \approx H$ and $\delta_1 = \delta_2, \beta_1 = \beta_2$ for all moderate functions $e$. In particular, $\delta_i = \beta_i = p$, $i = 1, 2$, for $e(r) = r^p, p > 0$ and for $a_p(t) = EX_t^{*p}$, $\lim_{t \to 0} a_p(t)^{-\alpha} X_t^* = 0$ or $\infty$ a.s. according as $\alpha < 1/p$ or $\alpha > 1/p$. The order of $a_p$ is known in special cases. Next we relate $Ee(|X_{T_r}|)$ to the moments of $T_r$. Pruitt [6] showed that $c_1 \leq Ey_{T_r}(r) = ET_r h(r) \leq c_2$ for Lévy processes. We extend the result to the present additive process setting. Define for $r > 0, n(r) = \inf\{t > 0 : y_t(r) > m\}$ where $m = (2\pi_d)^{-1}$.

LEMMA 5.1. *Let $X_t$ be a process in $\mathbb{R}^d$ with additive components. Then for any nondecreasing right-continuous function $\psi$ with $\psi(0) = 0$,*

(5.1) $$E\psi(T_r) \geq 2^{-1}\psi \circ n(r), \qquad r > 0.$$

PROOF. By the first bound in (2.6) and the definition of $n(r)$, $P(X_t^* \geq r) \leq \pi_d y_t(r) \leq \pi_d m = 1/2$ for $t \in [0, n(r)]$. Thus, $E\psi(T_r) = \int_0^\infty P(T_r \geq t)\psi(dt) = \int_0^\infty P(X_t^* \leq r)\psi(dt) \geq \int_0^{n(r)} P(X_t^* \leq r)\psi(dt) \geq 2^{-1} \int_0^{n(r)} \psi(dt) = 2^{-1}\psi \circ n(r)$. □

LEMMA 5.2. *Let $X_t$ be any process in $\mathbb{R}^d$ with additive components and let $e$ be any moderate function. Then there are two universal constants $c_1, c_2 \in (0, \infty)$ depending on $e$ and $d$ only such that*

(5.2) $$c_1 \leq Ee(y_{T_r}(r)) \leq c_2, \qquad r > 0.$$

PROOF. Let $\psi(t) = e(y_t(r))$ with $r$ fixed. By Lemma 5.1, $Ee(y_{T_r}(r)) = E\psi(T_r) \geq 2^{-1}\psi \circ n(r) = 2^{-1}e(m)$ since $y_{n(r)}(r) = m$. Since $\psi(t)/\psi(n(r)) \leq \rho(y_t(r)/y_{n(r)}(r))^\sigma$ for $t \geq n(r)$, that is, $\psi(t)/e(m) \leq \rho(y_t(r)/m)^\sigma$, $y_t(r)^{k/2} \geq \rho^{-(1+\delta)}m^{k/2}(\psi(t)/\ e(m))^{1+\delta}$ where $\delta = k/2\sigma - 1$ for some integer $k > 2\sigma$. The second bound in (2.6) yields

$$Ee(y_{T_r}(r)) = E\psi(T_r) = \int_0^{n(r)} P(X_t^* \leq r)\psi(dt) + \int_{n(r)}^\infty P(X_t^* \leq r)\psi(dt)$$

$$\leq e(m) + C' \int_{n(r)}^\infty (e(m)/\psi(t))^{1+\delta}\psi(dt)$$

$$= (1 + C'\delta^{-1})e(m) - C'\delta^{-1}e(m)^{1+\delta}\psi(\infty)^{-\delta} \leq (1 + C'\delta^{-1})e(m)$$

where $C' = A_k(d)\rho^{1+\delta}m^{-k/2}$. □

PROPOSITION 5.3. *Under the assumptions of Lemma 5.1, if there exists a nondecreasing continuous function $J$ with $J(0) = 0, J(t) > 0, t > 0, J(\infty) = \infty$ such that $y_{J(t_2)}(r)/y_{J(t_1)}(r) \geq \rho(t_2/t_1)^\sigma$ whenever $0 < t_1 < t_2$ for two constants $\rho, \sigma \in (0, \infty)$ not depending on $t_1, t_2, r$ and $\psi \circ J$ is moderate, then there is a constant $c \in (0, \infty)$ depending only on $\rho, \sigma, \psi \circ J$ and $d$ such that*

(5.3) $$E\psi(T_r) \leq c\psi \circ n(r), \qquad r > 0.$$

PROOF. Denote by $J^{-1}$ the inverse of $J$. Note that $J(J^{-1}(t)) = t$. Since $\psi \circ J$ is moderate, $\psi \circ J(t_2)/\psi \circ J(t_1) \leq \gamma(t_2/t_1)^\theta, 0 < t_1 < t_2$ with two constants $\gamma, \theta \in (0, \infty)$. As always,

$$E\psi(T_r) = \int_0^{J^{-1}(n(r))} P(X_{J(t)}^* \leq r)\psi \circ J(dt) + \int_{J^{-1}(n(r))}^\infty P(X_{J(t)}^* \leq r)\psi \circ J(dt)$$

$$\leq \psi \circ n(r) + \int_{J^{-1}(n(r))}^\infty P(X_{J(t)}^* \leq r)\psi \circ J(dt).$$

Let $\delta = k\sigma/2\theta - 1$ for an integer $k > 2\theta/\sigma$. The assumptions on $y_{J(t)}(r)$ and $\psi \circ J$ imply that $y_{J(t)}^{k/2}(r) \geq c_1(\psi \circ J(t)/\psi \circ n(r))^{1+\delta}$, where $c_1 = (\rho m\gamma^{-\sigma/\theta})^{k/2}$, for $t \geq J^{-1}(n(r))$. Thus, by the second bound in (2.6), $\int_{J^{-1}(n(r))}^\infty P(X_{J(t)}^* \leq r)\psi \circ J(dt) \leq c_2\delta^{-1}\psi \circ n(r)$, where $c_2 = A_k(d)\gamma^{1+\delta}(\rho m)^{-k/2}$. □

PROPOSITION 5.4. *Let $X_t$ be any additive process in $\mathbb{R}^d$ and let $e$ be any moderate function. Then*

(5.4) $$Ee(|X_{T_r}|) \approx e(r) + \int_r^\infty EG_{T_r}(\lambda)e(d\lambda), \qquad r > 0,$$

*where the constants in $\approx$ depend on $e$ and $d$ only.*

Proof. First we show that

(5.5) $$P(|\triangle X_{T_r}| > \lambda + 2r) = EG_{T_r}(\lambda + 2r), \quad \lambda > 0.$$

We have

$$P(|\triangle X_{T_r}| > \lambda + 2r, t < T_r \leq t + dt)$$
$$= P(|\triangle X_s| > \lambda + 2r \text{ for some } s \in (t, t+dt], X_t^* \leq r, |X_s| > r)$$
$$= P(|\triangle X_s| > \lambda + 2r \text{ for some } s \in (t, t+dt], X_t^* \leq r)$$

(if $|X_{s-}| \leq r$ and $|\triangle X_s| > \lambda + 2r$, then $|X_s| > r$

by the definition of $\triangle X_s = X_s - X_{s-}$)

$$= P(|\triangle X_s| > \lambda + 2r \text{ for some } s \in (t, t+dt]) \cdot P(X_t^* \leq r)$$

($X_s - X_t$ is independent of $\mathcal{F}_t, s > t$)

$$= (1 - \exp\{-(G_{t+dt}(\lambda + 2r) - G_t(\lambda + 2r))\}) \cdot P(T_r \geq t)$$

(the probability that $X$ has at least one jump of size larger than

$\lambda + 2r$ on $(t, t+dt]$ equals

$1 - \exp\{-(G_{t+dt}(\lambda + 2r) - G_t(\lambda + 2r))\}$;

the quasi-left-continuity of $X$, especially $\triangle X_t = 0$ a.s.

implies that $P(X_t^* \leq r) = P(T_r \geq t))$

$$= (G_{t+dt}(\lambda + 2r) - G_t(\lambda + 2r)) \cdot P(T_r \geq t) \quad (1 - e^{-x} \sim x \text{ as } x \downarrow 0).$$

Taking integration yields (5.5). Observe that $|X_{T_r}| + r \geq |\triangle X_{T_r}|$. Hence, $(|\triangle X_{T_r}| > 2r + \lambda) \subset (|X_{T_r}| - r > \lambda) \subset (|\triangle X_{T_r}| > \lambda)$ for $r > 0, \lambda > 0$. Also note that $Ee(k|X_{T_r}|) \approx Ee(|X_{T_r}|)$ for $k \in (0, \infty)$ where the constants in $\approx$ depend only on $e$ and $k$.

*Lower bound:*

$$Ee(3|X_{T_r}|) = \int_0^\infty P(3|X_{T_r}| \geq \lambda)e(d\lambda) = e(3r) + I_1$$

[since $|X_{T_r}| \geq r, P(3|X_{T_r}| \geq \lambda) = 1$ if $\lambda \leq 3r$] where $I_1 = \int_{3r}^\infty P(|X_{T_r}| \geq \lambda/3)e(d\lambda)$

$$\geq \int_{3r}^\infty P(|X_{T_r}| > \lambda/3 + r)e(d\lambda) \geq \int_{3r}^\infty P(|\triangle X_{T_r}| > \lambda/3 + 2r)e(d\lambda)$$

$$= \int_{3r}^\infty EG_{T_r}(\lambda/3 + 2r)e(d\lambda) \quad [\text{by (5.5)}]$$

$$\geq \int_{3r}^\infty EG_{T_r}(\lambda)e(d\lambda) \quad (\lambda/3 + 2r \leq \lambda \text{ because } \lambda \geq 3r).$$

By Lemma 5.2, $Ey_{T_r}(r) \leq c'$. Thus, $\int_r^{3r} EG_{T_r}(\lambda)e(d\lambda) \leq c \int_r^{3r} Ey_{T_r}(\lambda)e(d\lambda) \leq 2cEy_{T_r}(r) \int_r^{3r} e(d\lambda) < 2cc'e(3r)$. It follows that $e(3r) + I_1 \geq c_1(e(3r) + \int_r^\infty EG_{T_r}(\lambda)e(d\lambda))$.

*Upper bound:*

$$Ee(4^{-1}|X_{T_r}|) = \int_0^r P(4^{-1}|X_{T_r}| \geq \lambda)e(d\lambda) + I_2 \leq e(r) + I_2$$

where $I_2 = \int_r^\infty P(|X_{T_r}| \geq 4\lambda)e(d\lambda)$

$$\leq \int_r^\infty P(|X_{T_r}| > \lambda + 3r)e(d\lambda) \leq \int_r^\infty P(|\triangle X_{T_r}| > \lambda + 2r)e(d\lambda)$$

$$= \int_r^\infty EG_{T_r}(\lambda + 2r)e(d\lambda) \qquad [\text{by } (5.5)]$$

$$< \int_r^\infty EG_{T_r}(\lambda)e(d\lambda). \qquad \square$$

We can define

$$h(r) = \left[e(r) + \int_r^1 EG_{T_r}(\lambda)e(d\lambda)\right]^{-1}, \qquad r \in (0,1)$$

now. $H(r) \leq c(r)^{-1}h(r)$ if and only if there exists some $e$ moderate such that

(5.6) $$\int_r^1 EG_{T_r}(\lambda)e(d\lambda) \leq c(r)^{-1}e(r), \qquad r \in (0,1),$$

which is also equivalent to $Ee(|X_{T_r}|) \leq c(r)^{-1}e(r)$. Note that (5.6) has trivial solutions because many slow functions are moderate. Since $\int_r^1 EG_{T_r}(\lambda)e(d\lambda) \leq 2cEy_{T_r}(r)\int_r^1 e(d\lambda) < 2cc'e(1)$, (5.6) holds for $e(r) = c_1(r), c(r) = (2cc'c_1(1))^{-1}c_1(r)$ where $c_1(r)$ is any slow moderate function $[(\log(1/r))^{-1}$, say], which leads to $\delta_1 = \delta_2 = \beta_1 = \beta_2 = 0$. We rephrase all the above in one theorem.

THEOREM 5.5. *Let $X_t$ be any additive process in $\mathbb{R}^d$ and let $e$ be any moderate function. Define $u = a_e, h, H$ as above. Then (1.1) holds. If $EG_{T_r}(\lambda) \leq c(r)^{-1}A(r)/A(\lambda), 0 < r \leq \lambda \leq 1$ for a moderate function $A$, then (5.6) holds with $e = A^{1/q}$, $q > 1$, and hence $\delta_1 = \delta_2, \beta_1 = \beta_2$.*

The latter part holds because $\int_r^1 EG_{T_r}(\lambda)e(d\lambda) \leq c(r)^{-1}e(r)^q \int_r^1 e(\lambda)^{-q} \times e(d\lambda) = c_1 c(r)^{-1}e(r)^q (e(r)^{1-q} - e(1)^{1-q})$.

If $X_t$ is not a step process initially, (2.5) implies that $y_b(r)^{-1}, (Ey_{T_b}(r))^{-1}$ with fixed constants $b \in (0,\infty)$ are comparable to moderate functions. Thus, we can take $e(r) = y_b(r)^{-1}$ or $e(r) = (Ey_{T_b}(r))^{-1}$. For $e(r) = y_b(r)^{-1}$, $\delta_1 = \delta, \beta_2 = \beta$ where $\delta, \beta$ are the indices in Theorem 4.1. If $EG_{T_r}(\lambda) \leq c(r)^{-1}Ey_{T_b}(\lambda)/Ey_{T_b}(r)$ [resp. $EG_{T_r}(\lambda) \leq c(r)^{-1}y_b(\lambda)/y_b(r)$], $0 < r \leq \lambda \leq 1$ for some $b \in (0,\infty)$, then $\delta_1 = \delta_2, \beta_1 = \beta_2$ with $e(r) = (Ey_{T_b}(r))^{-1/q}$ [resp. $e(r) = y_b(r)^{-1/q}$], $q > 1$. In the case of Lévy processes, $EG_{T_r}(\lambda) = ET_rG(\lambda) \leq ET_rh(\lambda) \approx$ $h(\lambda)/h(r)$. Assume that $X_t$ is not a compound Poisson process and let $e = h^{-1/q}, q > 1$. To obtain Pruitt's result, we only need to prove that $a_e(t) \approx t^{1/q}$ where $e = \mathsf{k}^{1/q}$. Yang [8] showed that $a_e(t) \approx e \circ \mathsf{k}^{-1}(t) + t \int_{\mathsf{k}^{-1}(t)}^{\infty} G(\lambda) e(d\lambda)$ for all moderate $e$. For $e = \mathsf{k}^{1/q}$, $t \int_{\mathsf{k}^{-1}(t)}^{\infty} G(\lambda) e(d\lambda) \leq t \int_{\mathsf{k}^{-1}(t)}^{\infty} h(\lambda) e(d\lambda) \approx t \int_{\mathsf{k}^{-1}(t)}^{\infty} \mathsf{k}(\lambda)^{-1} e(d\lambda) = ct^{1/q}$. Thus, $a_e(t) \approx t^{1/q}$. In the more general case that $y_t(r) \approx f(t)z(r)$, $Ef(T_r) \approx z(r)^{-1}$ by Lemma 5.2. Thus, $EG_{T_r}(\lambda) \leq c_1 z(\lambda)/z(r)$ with $z(r)^{-1}$ comparable to a moderate function.

Revisit the example of $X_t^3 = X_{f_1(t)}^1 + X_{f_2(t)}^2$ with $y_t(r) = f_1(t)h_1(r) + f_2(t)h_2(r)$ for $X_t^3$ in Section 3: By Lemma 5.2, $Ey_{T_r}(r) = Ef_1(T_r)h_1(r) + Ef_2(T_r)h_2(r) \approx c \in (0, \infty)$. $EG_{T_r}(\lambda) \leq c_1 y_b(\lambda)/y_b(r)$ if $Ef_1(T_r)h_2(r) + Ef_2(T_r)h_1(r) \leq C < \infty$ for all small $r$.

Let $X_t$ be any additive process in $\mathbb{R}^d$. Disintegrating $\nu$ into $\nu(ds, dx) = \kappa_s(dx)g(ds)$ with a Lévy kernel $\kappa_s(dx)$ and a function $g$ yields $Q_t(r) = \int_{|x|>r} \nu_t(dx) + r^{-2} \int_{|x| \leq r} |x|^2 \nu_t(dx) = \int (|x|/r)^2 \wedge 1 \nu_t(dx) = \int_0^t \int (|x|/r)^2 \wedge 1 \kappa_s(dx) g(ds)$. If $\kappa_s$ satisfies $c_1 \theta(s) \kappa \leq \kappa_s \leq c_2 \theta(s) \kappa$ for a Lévy measure $\kappa$ and a function $\theta, c_1, c_2 \in (0, \infty)$ being two constants, that is, $\nu(ds, dx) \approx \theta(s)\kappa(dx)g(ds)$, then $Q_t(r) \approx g_1(t)Q(r)$ where $Q(r) = \int (|x|/r)^2 \wedge 1 \kappa(dx)$ and $g_1(t) = \int_0^t \theta(s)g(ds)$. Lemma 5.2 implies that $EG_{T_r}(\lambda) \leq c_3 Q(\lambda)/Q(r)$. $Q(r)^{-1}$ is continuous satisfying $Q(r) \leq C^2 Q(Cr)$ for $C > 1$ and is moderate if $\kappa$ is not a finite measure.

Take $e(r) = r^p, p > 0$ in Theorem 5.5. Then $\delta_1 = \beta_2 = p$. Note that $Ey_{T_r}(1) \geq 3^{-1} r^2 Ey_{T_r}(r) \approx r^2$ for $r \in (0, 1)$. If $Ey_{T_r}(1) \leq cr^2$, that is, $Ey_{T_r}(1) \approx r^2$, $\delta_1 = \delta_2 = \beta_1 = \beta_2 = p$ for all $p \in (0, 2)$. If $Ey_{T_r}(1) \leq cr^q, q \in (0, 2)$, $\delta_2 \geq q$ for all $p > 2$. The above can be derived from the crude estimate $\int_r^1 EG_{T_r}(\lambda)e(d\lambda) \leq c_1 \int_r^1 Ey_{T_r}(\lambda)e(d\lambda) \leq c_2 Ey_{T_r}(1) \int_r^1 \lambda^{-2} e(d\lambda)$.

$EG_{T_r}(\lambda)$ can be replaced by a function of $y_t(r)$. We have $EG_{T_r}(\lambda) = \int_0^\infty P(T_r \geq t) G(\lambda)(dt)$ while $P(T_r \geq t) = P(X_t^* \leq r) \leq c y_t(r)^{-k/2}$. Let $\nu(ds, dx) = \kappa_s(dx)g(ds)$. If $g$ is absolutely continuous, $G(\lambda)(dt) = g'(t)\tilde{G}_t(\lambda)\, dt$ where $\tilde{G}_t(\lambda) = \int_{|x|>\lambda} \kappa_t(dx)$. Hence $\int_r^1 EG_{T_r}(\lambda)e(d\lambda) = \int_0^\infty P(T_r \geq t) g'(t) \times (\int_r^1 \tilde{G}_t(\lambda)e(\lambda))\, dt$. But there is only minuscule gain in information on the structure of $e$ in (5.6) with that change. Constructing moderate functions $e$ satisfying $Ee(|X_{T_r}|) \leq c(r)^{-1} e(r)$ for small $r$, equivalently (5.6), remains open despite its enormous applications. [Obviously $Ee(|X_{T_r}|) \geq e(r)$.] A discrete construction method seems to be needed to tackle the problem. The case that $X_t$ is a step process initially can be shrugged off since $\inf_{1/2 \geq r > 0} \int_r^1 EG_{T_r}(\lambda)e(d\lambda) > 0$. One can also solve $\int_r^1 EG_{T_r}(\lambda)e(d\lambda) \leq e(r)$ or simply the integral-differential equation $\int_r^1 EG_{T_r}(\lambda)e'(\lambda)\, d\lambda = e(r)$ for a moderate function $e \in C^1$. Here we have an unbounded kernel $EG_{T_r}(\lambda)$. Since $Ee(|X_{T_r}|) \leq C_1 + C_2 e(r)$ implies $Ee(|X_{T_r}|) \leq C_3 e(r)$ for large $r$, the situation of (5.6) for large $r$ (with $\int_r^1$ replaced by $\int_r^\infty$) is slightly better where some results are available under the moment conditions.

Separately, given an additive process $X_t$, can we find the fully decomposed bounds

$$(5.7) \quad P(X_t^* \geq r) \leq c(r)^{-1} f(t) h(r), \qquad P(X_t^* \leq r) \leq (f(t) h(r))^{-1}$$

for some functions $f, h$ with a slow function $c(r)$? This is another important question yet to be answered. Equation (1.3) follows from (5.7) immediately. Equation (5.6) implies (5.7). Lemma 2.1 gives no information about (5.7) although $y_t(r) \approx f(t) z(r)$ holds in individual cases.

Schilling [7] deals with a class of Feller processes whose generators have the Lévy–Khintchine representation similar to the one used in additive processes with $B_t, Q_t, \nu_t, \Psi_t(\lambda)$ replaced by $B_x, Q_x, \nu_x, \Psi_x(\lambda)$, respectively, where $x = X_0$ and is mainly about the result that $P(X_t^* \geq r) \leq c_1 t H(r), P(X_t^* \leq r) \leq c_2 (th(r))^{-1}$ at $X_0 = x$ with $H, h$ defined in terms of $\Psi_x(\lambda)$ but the second bound requires the sector condition $|\mathrm{Im}\Psi_x(\lambda)/\mathrm{Re}\Psi_x(\lambda)| \leq c_0$, which probably can be removed if $h$ is defined by the characteristics not by the exponent. Schilling listed four possible cases of $\Psi_x(\lambda)$ ([7], Example 5.5 (4), (a)–(d), page 598) in which (1.3) holds. These cases are made with the common assumption that $\Psi_x(\lambda)$ can be decomposed into two elements, a Lévy exponent $\Psi(\lambda)$, and independently a function of $x$.

**6. Problem (1.3).** In the case of Lévy processes, $\Psi_t(\lambda) = t\Psi(\lambda)$ and we know the order of the function in the law of the iterated logarithm. We take $u(t) = t$ in (1.1) and have $\delta_2 = \delta_1, \beta_1 = \beta_2$. But for a general additive process, as changes occur at any given point in time, $\Psi_t(\lambda)$ shows no signs of the function $u$ in (1.1). Since $y_t(b)$ and $Ee(X_t^*)$ do not miss any of these instantaneous changes as $t \to 0$, $y_t(b)$ and $Ee(X_t^*)$ are two of the most likely benchmark functions for $u$. For $X_t$ of class $\mathcal{I}$, $u(t) = \int_0^t h_s(b)^{-1} y(b)(ds)$ with $\int_0^t y(b)(ds) = y_t(b)$. We discuss below only the $\delta$-indices. $\beta$-indices follow suit.

Let $v$ be the function determined by equation $y_{v(t)}(b) = \phi(t)$ with $b \in (0, \infty)$ and $\phi$ quasiconvex. Which combination of $b$ and $\phi$ is the best? What we have is either $\delta_2 < \delta_1$ or $\delta_2 = \delta_1$. $\phi$ is probably more important than anything else but we do not know what kind of $\phi$ can make a change from $\delta_2 < \delta_1$ to $\delta_2 = \delta_1$ since we are unable to calculate $\delta_1, \delta_2$ in general. If $\delta_2 < \delta_1$, $\delta_1 - \delta_2$, big or small, makes no difference as long as $\delta_2 > 0$. Here, for instance, $u(t) = y_t(b)^{1/p}$ for $\phi(t) = t^p, p > 0$. In this case, if $\delta_2 < \delta_1$ (resp. $\delta_2 = \delta_1$) for *some* $p$, then $\delta_2 < \delta_1$ (resp. $\delta_2 = \delta_1$) for *all* $p$. We consider $\phi$ as *acceptable* if $\delta_2 > 0$. Again, this is a technical matter. It is not easy to show that $\delta_2 > 0$. We have $\delta_\varepsilon \uparrow \delta^* \leq \delta_2$ [see (i), proof of Theorem 3.2] and $\delta_\varepsilon > 0$ for some $\varepsilon$ if and only if there exist $\alpha_1, \alpha_2 \in (0, \infty)$ such that $y_{v(r)}(r^{\alpha_1}) \geq r^{-\alpha_2}, r \in (0, r_1)$. One can replace $\delta_2$ by $\delta^*$ in Theorem 3.1. $\delta_\varepsilon$ and $\delta_1$ look similar without $\Sigma$-sequences and infinite sums in their definitions.

It is almost impossible to work on *infinitely often* events directly. We do not know how to do the following: Given a nondecreasing function $v$, decide whether $\overline{\delta} < \underline{\delta}$ or $\overline{\delta} = \underline{\delta}$. If any, find an example where $\overline{\delta} < \underline{\delta}$ for all nondecreasing functions $v$. (The case $\underline{\delta} = 0$ is excluded.) Other possible cases include $\delta_2 < \overline{\delta} = \underline{\delta}$, $\delta_2 = \overline{\delta} < \underline{\delta}$, $\delta_2 = \dot{\delta}_2 < \overline{\delta} = \underline{\delta}$, and so on. The only tool available is still the classical Borel–Cantelli–Fatou argument. There are some possible ways to construct $v$ such that $\overline{\delta} = \underline{\delta}$.

(a) At this juncture, it is unclear whether or not (1.3) has anything to do with $y_t(r)$. What is clear, though, is that (1.3) holds if $y_t(r)$ has a good structure. Consider a function $g_\varepsilon(r) > 0$ with the two properties: For each $\varepsilon$ there is a $\Sigma$-sequence $\sigma_n \downarrow 0$ such that $\sum g_\varepsilon(\sigma_n) < \infty$; $\lim_{\varepsilon \to 0} g_\varepsilon(r) = 1$ for every $r$; for example, $g_\varepsilon(r) = r^{c\varepsilon}, c > 0$. Define $\delta_\varepsilon = \inf\{\eta > 0 : \liminf_{r \to 0} g_\varepsilon(r) y_{v(r)}(r^{1/\eta}) = 0\}$, where $v$ is a nondecreasing function. The first property of $g_\varepsilon(r)$ implies that $\delta_\varepsilon \leq \delta_2$. With the second property, to get $\delta_1 = \delta_2$, we need to construct the function $v$ and the function $g_\varepsilon(r)$ from $y_t(r)$ such that $\lim_{n \to \infty} \delta_{\varepsilon_n} = \inf\{\eta > 0 : \liminf_{r \to 0} \lim_{n \to \infty} g_{\varepsilon_n}(r) y_{v(r)}(r^{1/\eta}) = 0\}$ for some sequence $\varepsilon_n \downarrow 0$. Quasiconvex condition (3.1) makes one such case.

(b) Like the fully decomposed bounds in (5.7), the following types of bounds also lead to $\overline{\delta} = \underline{\delta}$. Type I: $P(X_t^* \geq r) \leq u(t) c(r)^{-1} w_t(r), P(X_t^* \leq r) \leq w_t(r)^{-1}$, where $u$ is a nondecreasing continuous function with $u(0) = 0$. Let $v$ be the inverse of $u$ and define $\delta = \inf\{\eta > 0 : \liminf_{r \to 0} r^\varepsilon w_{v(r)}(r^{1/\eta}) = 0\}, \varepsilon \in (0, 1)$. Type II: $P(X_t^* \leq r) \leq g(r) c(t)^{-1} w_t(r)^{-1}, P(X_t^* \geq r) \leq w_t(r)$, where $g$ is a function satisfying $\sum g(\sigma_n) < \infty$ for some $\Sigma$-sequence $\sigma_n \downarrow 0$. Let $v$ be any nondecreasing function satisfying $w_{v(t_n)}(t_n^{1/\kappa}) \to 0$ for some sequence $t_n \downarrow 0$ and $\kappa \in (0, \infty)$ and define $\delta = \inf\{\eta > 0 : \liminf_{r \to 0} w_{v(r)}(r^{1/\eta}) = 0\}$. $c(r)$ stands for slow functions as always. If the bounds of the above types are obtainable, the reader can check that $\liminf_{t \to 0} t^{-1/\eta} X_{v(t)}^* = 0$ or $\infty$ a.s. according as $\eta > \delta$ or $\eta < \delta$. The above holds true for any appropriate process in $\mathbb{R}^d$. Normally, the result of de la Peña and Eisenbaum is relatively sharp in the case of additive processes. We wonder if it can be used to obtain the bounds in Type I with $u(t) = a_e(t)$ where $e$ is not a slow function. To do so, it is necessary to have a new mechanism of getting around the point where Markov's inequality enters. The proposed bounds can be called skew subdecomposable. We have been unable to find a convincing number of their examples. Even in the case of Lévy processes, it is still a question that such bounds can exist.

(c) Again the following holds true for any process $X_t$ in $\mathbb{R}^d$ continuous in probability with $X_0 = 0$. Both $P(X_t^* \geq r)$ and $P(X_t^* \leq r)$ diverge as $t, r$ move in tandem toward 0. This is a useful fact. The more we know about $P(X_t^* \geq r)$ and $P(X_t^* \leq r)$, the more likely to construct $v$ such that $\overline{\delta} = \underline{\delta}$. Fix $\delta \in (0, \infty)$. (i) We can construct two sequences $\bar{\sigma}_n \downarrow 0, v_n \downarrow 0$ such that $\sum P(X^*_{v_n} \le \bar{\sigma}_n^{1/\delta}) < \infty$ and $P(X^*_{v_{n+1}} \ge \bar{\sigma}_n^{1/\delta}) \to 0$. (ii) For any sequence $v_n \downarrow 0$ selected, we can construct a $\Sigma$-sequence $\bar{\sigma}_n \downarrow 0$ such that $\sum P(X^*_{v_n} \le \bar{\sigma}_n^{1/\delta}) < \infty$. (This time, $v_n \downarrow 0$ is the one fixed first.) Although we have assumed in (ii) that $P(X^*_t \le \varepsilon) \downarrow 0$ as $\varepsilon \downarrow 0$, there is no loss of generality because the assumption fails only when $X_t = 0$ a.s. for a length of time at first, in which case (1.3) is uninteresting. Part (ii) has also proved that there is always a strictly increasing continuous function $v$ with $v(0) = 0$ such that $\dot{\delta}_2 \ge \delta$ although there is no guarantee for $\underline{\delta} < \infty$. If we can in (i) make $\{\bar{\sigma}_n\}$ a $\Sigma$-sequence or in (ii) make $\{v_n, \bar{\sigma}_n\}$ contain subsequences $\{v_{n_k}, \bar{\sigma}_{n_k}\}$ such that $P(X^*_{v_{n_k+1}} \ge \bar{\sigma}_{n_k}^{1/\delta}) \to 0$ as $k \uparrow \infty$, then $\liminf_{t \to 0} t^{-1/\eta} X^*_{v(t)} = 0$ or $\infty$ a.s. according as $\eta > \delta$ or $\eta < \delta$ for any nondecreasing function $v$ with $v(\bar{\sigma}_n) = v_n$. The first result follows from the fact that $X^*_{v_{n_k+1}}/\bar{\sigma}_{n_k+1}^{1/\eta} = (X^*_{v_{n_k+1}}/\bar{\sigma}_{n_k}^{1/\delta}) \cdot ((\frac{\bar{\sigma}_{n_k}}{\bar{\sigma}_{n_k+1}})^{1/\eta} \bar{\sigma}_{n_k}^{1/\delta - 1/\eta})$ ($\eta > \delta$) where the first factor is bounded by 1 for infinitely many $k$ and the second tends to 0 since $\{\bar{\sigma}_n\}$ is a $\Sigma$-sequence. It is also clear from here that $\underline{\delta} \le \delta$. The second result is obvious because $\eta < \delta \le \dot{\delta}_2$. Hence $\dot{\delta}_2 = \overline{\delta} = \underline{\delta} = \delta$.

As this paper draws to a close, we say a few words about changing measures. It is a technique used to identify more processes for which (1.3) holds. Let $(X_t; B_t, \nu_t)$ be an additive process in $\mathbb{R}^d$ and $(X'_t; B'_t, \nu'_t)$ another. $X_t$ (resp. $X'_t$) induces a probability measure $P_t$ (resp. $P'_t$) on the canonical space. $P'_t \ll P_t$, that is, $P'_t$ is absolutely continuous with respect to $P_t$, if and only if there exists a Borel function $f:[0,t] \times \mathbb{R}^d \to \mathbb{R}_+$ satisfying $\int_{[0,t] \times \mathbb{R}^d} (1 - \sqrt{f})^2 d\nu < \infty$ such that $\nu'(ds, dx) = f(s,x)\nu(ds, dx)$ and $B'_s = B_s + \int_0^s \int_{|x| \le 1} (f(\tau, x) - 1) x \nu(d\tau, dx), s \in (0, t]$. This result can be found in [3], Chapter IV, and [2], Chapter XIV. We did not take $Q_t$ into account in order to simplify the formula. If $f > 0$, $P_t \ll P'_t$ as well with $f$ replaced by $1/f$, that is, $P'_t \sim P_t$. Here, for instance, $\int_{[0,t] \times \mathbb{R}^d} (1 - \sqrt{f})^2 d\nu < \infty$ for $f = e^{\phi}$ with $\phi:[0,t] \times \mathbb{R}^d \to [-C, C]$ satisfying $|\phi(s,x)| \approx |x|$ near $x = 0$. Normally, one relates $B_t, \nu_t$ to $\phi$. If $P'_t \ll P_t$, $X'_t$ has the sample-path behavior of $X_t$ up to time $t$ at least even though $X'_t$ and $X_t$ are totally different in law. If $P'_t \ll P_t$, $\gamma'_0(t) = \gamma_0(t)$ and monotonicity of the sample-path is preserved, that is, the drift remains nondecreasing. There are also some other invariant properties that we omit. Therefore, $P'_t \ll P_t$ cannot occur unless and until the two additive processes in consideration belong to the same category. Let $y'_t(r)$ denote the function in (2.4) for $X'_t$. Componentwise, $y'_t(r) = G'_t(r) + K'_t(r) + M'^*_t(r)$ with $G'_t(r) + K'_t(r) = \int_{[0,t] \times \mathbb{R}} (x/r)^2 \wedge 1 f(s,x) \nu(ds, dx)$, $M'_t(r) = |m'_t(r)|$ where $m'_t(r) = m_t(r) + r^{-1} \int_0^t \int_{|x| \le r} (f(\tau, x) - 1) x \nu(d\tau, dx)$ and $M_t(r) = |m_t(r)|$. Usually, there are no similarities between $y'_t(r)$ and $y_t(r)$. If $X_t$ satisfies (1.3), so does $X'_t$ while $v_i, v_s, \delta, \beta$ need no change, but $y'_t(r)$ may not satisfy any condition for (1.3) so far developed in this paper. [In other words, the shape of the function in (2.4) is only one piece of the puzzle.] Similarly, if $X_t$ is a semimartingale additive process of class $\mathcal{I}$, then (1.3) holds for $X'_t$ where $\nu'(ds,dx) = \kappa'_s(dx)u(ds)$, $\kappa'_s(dx) = f(s,x)\kappa_s(dx)$, $B'_t = \int_0^t b'_s u(ds)$, $b'_s = b_s + \int_{|x|\le 1}(f(s,x)-1)x\kappa_s(dx)$, $\tilde{m}'_s(r) = \tilde{m}_s(r) + r^{-1}\int_{|x|\le r}(f(s,x)-1)x\kappa_s(dx)$. Conditions (i) and (ii) may fail to hold for $X'_t$. Still, if $X_t = \bar{X}_{g(t)}$ where $(\bar{X}; B, \nu)$ is a Lévy process in $\mathbb{R}^d$ and $g$ is a nondecreasing continuous function with $g(0) = 0$, then (1.3) holds for $X'_t$ (as well as for $X_t$) with $v_i(t) = v_s(t) = g^{-1}(t) = \inf\{s > 0 : g(s) > t\}$ since $y_t(r) = g(t)h(r)$. $X'_t$ having characteristics $\nu'(ds,dx) = f(s,x)g(ds)\nu(dx)$, $B'_t = g(t)B + \int_0^t \int_{|x|\le 1}(f(s,x)-1)x\nu(dx)g(ds)$ is a genuine nonhomogeneous process. Perhaps $X'_t$ satisfying (1.3) can only be recognized through changing measures. The same occurs for Lévy processes. Let $(X_t, \nu)$ be a strictly $\alpha$-stable process in $\mathbb{R}^d$. For $d = 1$, $\nu$ takes the form $\nu(dx) = s(x)\,dx$ with $B$ determined by $\nu$ for $\alpha \ne 1$. Any Borel function $\theta > 0$ satisfying $\int(\sqrt{\theta(x)} - \sqrt{s(x)})^2\,dx < \infty$ defines a Lévy process $X'_t$ with $\nu'(dx) = f(x)\nu(dx) = \theta(x)\,dx, B' = B + \int_{|x|\le 1}(\theta(x) - s(x))x\,dx$, where $f(x) = \theta(x)/s(x)$. $X'_t$ nowhere near stable has exactly the same local sample-path behavior as $X_t$ does while many fine results on the sample paths hold for stable processes. $X'_t$ becomes nonhomogeneous if $f(x)$ is replaced by $f(s,x)$.

The measure change formula also opens up a way to prove (1.3) for $X_t$, that is, to find a function $f > 0$ in the formula such that one of the conditions for (1.3) can hold for $y'_t(r)$. Take, for example, a symmetric additive process $X_t$ in $\mathbb{R}$ with $y_t(r) = \int_{[0,t]\times\mathbb{R}}(x/r)^2 \wedge 1\nu(ds,dx) = \int_0^t \int (x/r)^2 \wedge 1\kappa_s(dx)u(ds)$. Let $f = e^\phi$, $\phi(s,x) = \theta(s)(|x|\wedge 1)$ where $\theta(s) = \theta(s, \nu_s, \kappa_s, u(s))$ is any bounded Borel function. Then $X'_t$ is also symmetric with $y'_t(r) = \int_0^t k(s,r)u(ds)$ where $k(s,r) = \int (x/r)^2 \wedge 1 e^{\theta(s)(|x|\wedge 1)}\kappa_s(dx)$. There are several options in constructing $\theta$. Here, for instance, if $k(s,r) \approx l(s)z(r)$, $y'_t(r) \approx g(t)z(r)$ with $g(t) = \int_0^t l(s)u(ds)$; if $k(s,r)$ is nondecreasing in $s$, $X'_t$ is of class $\mathcal{I}$; if $k(t,r)/k(t,r_1) \ge y'_t(r)/y'_t(r_1)$ for $r_1 < r$ with $u$ differentiable, $y'_t(r)/y'_t(r_1)$ is nondecreasing in $t$ for $r_1 < r$, so (i) of Proposition 4.2 holds.

**Acknowledgment.** The author would like to thank the referee for his careful reading of the paper and for valuable suggestions from which necessary expository improvements have been made.

P.O. BOX 647
JACKSON HEIGHTS, NEW YORK 11372
USA
E-MAIL: myang1968@yahoo.com